\documentclass[10pt,leqno]{article}

\usepackage{amssymb}
\usepackage{amsmath}
\usepackage{amsfonts}
\usepackage{amscd}
\usepackage{amsthm}
\usepackage{graphicx}
\usepackage{enumerate}
\usepackage{empheq}
\usepackage{xcolor}

\newtheorem{theorem}{Theorem}[section]

\newtheorem{Remark}{Remark}[section]
\newtheorem{lemma}{Lemma}[section]
\newtheorem{example}{Example}[section]

\newcommand{\se}{\setcounter{equation}{0}}

\newcommand{\R}{\mathbb{R}}

\newcommand{\bu}{\textbf{u}}
\newcommand{\bU}{\textbf{U}}
\newcommand{\bW}{\textbf{W}}
\newcommand{\bX}{\textbf{X}}

\newcommand{\bV}{\textbf{V}}
\newcommand{\bS}{\textbf{S}}
\newcommand{\bR}{\textbf{R}}
\newcommand{\bP}{\textbf{P}}
\newcommand{\bH}{\textbf{H}}
\newcommand{\bL}{\textbf{L}}
\newcommand{\bv}{\mathbf{v}}
\newcommand{\bw}{\textbf{w}}
\newcommand{\de}{\textbf{e}}
\newcommand{\bn}{\textbf{n}}
\newcommand{\bz}{\textbf{z}}
\newcommand{\bq}{\textbf{q}}
\newcommand{\bJ}{\textbf{J}}
\newcommand{\brf}{\mathbf{f}}
\newcommand{\al}{\alpha}
\newcommand{\be}{\mbox{\boldmath$\eta$}}
\newcommand{\bphi}{\mbox{\boldmath$\phi$}}
\newcommand{\bxi}{\mbox{\boldmath$\xi$}}
\newcommand{\bs}{\mbox{\boldmath$\psi$}}
\newcommand{\bt}{\mbox{\boldmath$\theta$}}
\newcommand{\bzeta}{\mbox{\boldmath$\zeta$}}
\newcommand{\bchi}{\mbox{\boldmath$\chi$}}

\newcommand{\vertiii}[1]{{\left\vert\kern-0.25ex\left\vert\kern-0.25ex\left\vert #1 
		\right\vert\kern-0.25ex\right\vert\kern-0.25ex\right\vert}}

\begin{document}
	\author{Saumya Bajpai\footnotemark[1]\thanks{School of Mathematics and Computer Science, Indian Institute of Technology Goa, Ponda, Goa-403401, India. Email: saumya@iitgoa.ac.in
},~ Deepjyoti Goswami\footnotemark[2]\thanks{Department of Mathematical Sciences, Tezpur University, Tezpur, Sonitpur, Assam-784028, India. Email: deepjyoti@tezu.ernet.in, kallol@tezu.ernet.in},~ and Kallol Ray\footnotemark[2]}
	\title{
		Optimal Error Estimates of a Discontinuous Galerkin Method for the Navier-Stokes Equations}
	\date{}
	\maketitle

	\begin{abstract}
		In this paper, we apply discontinuous finite element Galerkin method to the time-dependent $2D$ incompressible Navier-Stokes model. We derive optimal error estimates in $L^\infty(\bL^2)$-norm for the velocity and in $L^\infty(L^2)$-norm for the pressure  with the initial data $\bu_0\in \bH_0^1\cap \bH^2$ and the source function $\brf$ in $L^\infty(\bL^2)$ space. These estimates are established with the help of  a new $L^2$-projection and modified Stokes operator on appropriate broken Sobolev space and with standard parabolic or elliptic duality arguments. Estimates are shown to be uniform under the smallness assumption on data. Then, a completely discrete scheme based on the backward Euler method is analyzed, and fully discrete error estimates are derived. We would like to highlight here that the estiablished semi-discrete error estimates related to the $L^\infty(\bL^2)$-norm of velocity  and $L^\infty(L^2)$-norm of pressure are optimal and sharper than those derived in the earlier articles. Finally, numerical examples validate our theoretical findings.

	\end{abstract}

	\vspace{1em} 
	
	\section{Introduction}\label{s1}
	\se
	In this paper, we consider the  following time-dependent Navier-Stokes equations for incompressible fluid flow:
	\begin{eqnarray}
		\textbf{u}_t-\mu\Delta \textbf{u}+ \textbf{u}\cdot \nabla \textbf{u}+\nabla p &= & \textbf{f} \quad \text{in $\Omega$ for $0<t\leq T$,}\label{8.1} \\ 
		\nabla\cdot \textbf{u} &=&  0 \quad \text{in $\Omega$ for $0<t\leq T$,}\label{8.2}\\ 
		\textbf{u} &=& \textbf{u}_0 \quad \text{in $\Omega$ for $t=0$,}\label{8.3}\\
		\textbf{u} &=  & 0 \quad \text{on $\partial\Omega$ for $0<t\leq T$},\label{8.4}
	\end{eqnarray}
	where $\bu$ is the fluid velocity, $p$ the pressure, $\brf$ the external force, $\mu>0$ the kinematic viscosity, and $\Omega\subset\R^2$ a bounded, simply connected domain with polygonal boundary $\partial\Omega$. And $T>0$ is the final time. We also impose the usual normalization condition on the pressure, namely, that $\displaystyle{\int_\Omega} p\, d\Omega=0$. 
	Throughout the article, the boldface letters denote the vectors.
	
	There is a large amount of literature devoted to finite element and related methods for the Navier-Stokes equations. However, very few are directed to the analysis of discontinuous Galerkin (DG) methods. These methods have been evolved to compensate the continuous Galerkin methods which fail in case of higher-order approximation and for unstructured grids. DG methods are known to be flexible in handling local mesh adaptivity and non-uniform degrees of approximation of solutions with variable smoothness.
	Besides, they are elementwise conservative and easy to implement than the finite volume methods and the standard mixed finite element methods when high degree piecewise polynomial approximations are involved.
	Introduced in \cite{LR74,RH73}, DG methods have been attempted to the Euler and  Navier-Stokes equations (NSEs) as early as in \cite{BO1999,LC93}. However, the first rigorous analysis in this direction can be attributed to Girault {\it et al.} \cite{GRW005}, where a DG method have been formulated
	with nonoverlapping domain decomposition, for the steady incompressible Stokes and Navier Stokes system of equations, with approximations of order $k=1,2$ or $3$.  The authors have  established a uniform discrete inf-sup condition and then have demonstrated optimal energy estimates for the velocity and $L^2$-estimates for the pressure.  A follow-up of the work \cite{RG06} has an improved inf-sup condition and has discussed several numerical schemes and numerical convergence rates.  The extension to time-dependent NSEs can be found in \cite{KR05}, where an error analysis of a  linear subgrid-scale eddy viscosity method  combined with discontinuous Galerkin approximations has been worked out. Optimal semi-discrete error estimates of the velocity and pressure have been derived with reasonable dependence on Reynolds number. Then two fully discrete schemes, which are first and second-order in time, respectively, have been analyzed and, optimal velocity error estimates have been established. In \cite{GRW05}, Girault {\it et al.} have analyzed a projection method, decoupling the velocity and the pressure, along with a discontinuous Galerkin method for the time-dependent incompressible NSEs. Optimal error estimates for the velocity and suboptimal for pressure have been presented. Several other notable works can be found in \cite{CW10,CKS09,PE10,JHYW18} and references therein that deal with DG methods for incompressible NSEs. And for work based on numerical schemes and numerical convergence rates, we refer to  \cite{FW11, LFAR19, PMB18, SFE07}.
	For steady and unsteady incompressible Navier-Stokes problem, the error analysis more or less stops at energy error estimate for the velocity. Whereas numerically optimal rate of convergence for velocity, both in energy and in  $\bL^2$-norms has been shown, to the best of our knowledge no analysis is available for optimal $L^\infty(\bL^2)$-norm error estimate for velocity, except in  \cite{GRW005}. In \cite{GRW005}, for the steady state Navier-Stokes problem, authors have given a sketch of how to go about the proof of $\bL^2$-norm error estimate of velocity for SIPG (symmetric interior penalty Galerkin) method. However, no such result exists for time dependent Navier-Stokes problem. We would like to point out here that for NIPG (non-symmetric interior penalty Galerkin) method, the $\bL^2$-norm error convergence rate will depend on the degree of polynomial approximation; it will provide optimal and sub-optimal results, when the degree of polynomial is odd and even, respectively \cite{GRW05,R08,RG06}. For SIPG, the error estimate will not depend on the polynomial degree, but will depend on the penalty parameter $\sigma_e$ (see {\bf Section} \ref{s2}) which  has to be sufficiently large \cite{R08}. We have mainly considered here the SIPG method, since NIPG case as mentioned is known to give sub-optimal estimates in general.  The goal of this paper is to establish  analytically, the optimal $L^\infty(\bL^2)$-norm error estimate for velocity and  $L^\infty(L^2)$-norm error estimate for pressure, which does not follow immediately just by following the steady case, outlined in \cite{GRW005}, and is more technical than we have expected. 
	
	
	In the DG literature, the error analysis revolves around the 
	operator $\bR_h$ (see {\bf Section} \ref{s3})
	which is used to obtain optimal error estimates of the velocity in the evergy norm and of the pressure in $L^2$-norm. Using duality argument the optimal $L^\infty(\bL^2)$-norm error estimate for the velocity of the steady Navier-Stokes problem can be obtained (see \cite{GRW005}). However this procedure fails in the case of unsteady NSEs and we feel that this is due to the lack of appropriate approximation operator and projection in the DG finite element set-up. Following  the work of Heywood and  Rannacher for NSEs in a non-conforming set-up \cite{HR82}, we have first constructed an $L^2$-projection $\bP_h$ onto an appropriate DG finite element space, with the help of the operator $\bR_h$ and have proved the requisite approximation properties for $\bP_h$. Next, we have defined an approximation operator $\bS_h$, which is a modified Stokes projection (see {\bf Section} \ref{s3} for more details), again in DG set-up, that is, for broken Sobolev space. This is crucial for obtaining optimal error estimate using duality argument, as has been the standard procedure for parabolic problem or time dependent Stokes problem (linear NSE) in conforming Galerkin methods. Optimal approximation estimates for $\bS_h$ have been derived with the help of the projection $\bP_h$. Although we have followed the ideas of \cite{HR82} but there are differences and technical difficulties since our formulations and finite element spaces are different. For example, special care is needed to handle the nonlinear convection term.  Armed with these approximation operators, we have achieved here for $t>0$, optimal $L^2$ error estimates for the velocity and the pressure. 
	%
	Results obtained here are sharper than the existing ones derived for the finite element discontinuous Galerkin method applied to (\ref{8.1})-(\ref{8.4}).

	Below, we summarize our major results obtained in this article:
	\begin{itemize}
		\item New $L^2$-projection $\bP_h$ and an approximation operator (modified Stokes projection) $\bS_h$ on a suitable DG finite element space are introduced, and their approximation properties are derived.
		\item Optimal $L^\infty(\bL^2)$ and $L^\infty(L^2)$-norms error estimates for semi-discrete discontinuous Galerkin approximations to the velocity and the pressure, respectively, are derived, considering  $\brf,\brf_t\in L^\infty(\bL^2)$ and $\bu_0\in \bH^2\cap \bJ_1$  (see {\bf Section} \ref{s2}, for the definition of $\bJ_1$).
		\item Under smallness condition on the given data, uniform in time semidiscrete optimal error estimate for velocity is established.
		\item  $L^2$ error estimates for fully discrete DG approximations to the velocity and the pressure, respectively, are derived, when a first order backward Euler time discretization scheme is applied.%
	\end{itemize}

	The outline of the article is as follows. Notations,  variational formulation, and basic assumptions are presented in {\bf Section} \ref{s2}. {\bf Section} \ref{s3} contains a semi-discrete DG method and the standard projection properties. {\bf Section} \ref{s4} is devoted to the optimal $L^\infty(\bL^2)$-norm error estimates of the velocity term,  based on new projection properties. Further, under smallness condition on the data, uniform in time $t>0$, error estimates are established.  And in {\bf Section} \ref{s5}, optimal error estimates for the pressure are derived. In {\bf Section} \ref{s6}, a fully discrete scheme based on the backward Euler method is formulated, and error estimates for the velocity and pressure are derived.  Further, numerical experiments are conducted and the results obtained are analyzed in {\bf Section} \ref{s7}.
	
	\section{Preliminaries and discontinuous variational formulation}\label{s2}
	\se
	In the rest of the paper, we denote by bold face letters the $\R^2$-valued
	function spaces such as $\bH_0^1 = (H_0^1(\Omega))^2$, $ \bL^2 = (L^2(\Omega))^2$ etc.
	%
	The norm in $W^{l,r}(\Omega)$ is denoted by $\|\cdot\|_{l,r,\Omega}$ and the seminorm by $|\cdot|_{l,r,\Omega}$. The $L^2$ inner-product is denoted by $(\cdot,\cdot)$.  For the Hilbert space $H^l(\Omega):=W^{l,2}(\Omega)$, the norm is denoted by $\|\cdot\|_{l,\Omega}$ or $\|\cdot\|_l$.
	The space $\bH^1_0$ is equipped with the norm
	\[ \|\nabla\bv\|= \big(\sum_{i,j=1}^{2}(\partial_j v_i, \partial_j
	v_i)\big)^{1/2}=\big(\sum_{i=1}^{2}(\nabla v_i, \nabla v_i)\big)^{1/2}. \]
	Let $H^m/\R$ be the quotient space with norm $\|\phi\|_{H^m/\R}=\inf_{c\in\R}\|\phi+c\|_m.$ For $m=0$, it is denoted by $L^2/\R$. 
	For any Banach space $X$, let $L^p(0, T; X)$ denote the space of measurable $X$-valued functions $\bphi$ on  $ (0,T) $ such that
	\[ \int_0^T \|\bphi (t)\|^p_X~dt <\infty~~\mbox {if}~1 \le p < \infty, ~~~\mbox{and}~~~
	{\displaystyle{ess \sup_{0<t<T}}} \|\bphi (t)\|_X <\infty~~\mbox {if}~p=\infty. \]
	The dual space of $H^m(\Omega)$, denoted by $H^{-m}(\Omega)$, is defined as the completion of 
	$C^\infty(\bar{\Omega})$ with respect to the norm 
	$$\|\phi\|_{-m}:=\sup \big\{\frac{(\phi,\psi)}{\|\psi\|_m}: \psi \in H^m(\Omega), \|\psi\|_m\neq 0 \big\}.$$
	%
	Further, we introduce some divergence free function spaces for our future use:
	\[ \bJ_1=\{\bw\in \bH_0^1(\Omega):\,\nabla\cdot \bw=0\}\quad \text{and}\quad \bJ=\{\bw\in \bL^2(\Omega):\,\nabla\cdot \bw=0,\, \bw\cdot \bn|_{\partial\Omega}=0~~\text{holds weakly}\},\] 
	where $\bn$ is the outward normal to the boundary $\partial\Omega$ and $\bw\cdot \bn|_{\partial\Omega}=0$ should be understood in the sense of trace in $\bH^{-1/2}(\partial\Omega)$. \\
	We now present below, a couple of assumptions on the given data and the domain.\\
	({\bf A1}). The initial velocity $\bu_0$ and the external force $\brf$ satisfy, for some positive constant $M_0$ and for time $T$ with $0<T<\infty$,
	
	$\bu_0\in \bH^2\cap \bJ_1,~~\brf,\brf_t\in L^\infty(0,T;\bL^2)$ with $\|\bu_0\|_2\leq M_0$ and $\displaystyle\sup_{0<t<T}\{\|\brf(\cdot,t)\|,\|\brf_t(\cdot,t)\|\}\leq M_0.$ \\
	({\bf A2}). For ${\bf g}\in \bL^2$, let the unique solutions $\bv\in \bJ_1,~q\in L^2/ \R$ of the steady state Stokes problem
	\begin{align*}
		-\Delta \bv+\nabla q={\bf g},\\
		\nabla\cdot \bv=0 ~\text{in}~\Omega,~~\bv|_{\partial\Omega}=0,
	\end{align*}   
	satisfy 
	\[\|\bv\|_2+\|q\|_{H^1/ \R}\leq C\|{\bf g}\|.\]
	We  next turn to the weak formulation of (\ref{8.1})-(\ref{8.4}) as follows: Find a pair $(\bu(t),p(t))\in \bH_0^1\times L^2/ \R$, $t>0$, such that
	\begin{align}
		(\bu_t,\bphi)+\mu(\nabla\bu,\nabla\bphi)+(\bu\cdot\nabla\bu,\bphi)-(p,\nabla\cdot\bphi)&=(\brf,\bphi)\quad \forall \bphi\in \bH_0^1,\label{w1}\\
		(\nabla\cdot\bu,q)&=0\quad \forall q\in  L^2/ \R,\label{w2}\\
		(\bu(0),\bphi)&=(\bu_0,\bphi)\quad \forall \bphi\in \bH_0^1.\label{w3}
	\end{align}
	And for DG formulation, we consider a regular family of triangulations $\mathcal{T}_h$ of $\bar{\Omega}$,
	consisting of triangles of maximum diameter $h$. 
	We also assume that the subdivision is regular. Let $h_T$ denotes the diameter of a triangle $T$ and $\rho_T$ the diameter of its inscribed circle.  Then, by regular, we mean that there exists a constant $\gamma>0$ such that
	\begin{equation}
		\frac{h_T}{\rho_T}\leq \gamma, \quad \forall T\in \mathcal{T}_h.\label{regular}
	\end{equation}
	We denote by $\Gamma_h$ the set of all edges of $\mathcal{T}_h$. With each edge $e$, we associate a unit normal vector $\bn_e$. If $e$ is on the boundary $\partial\Omega$, then $\bn_e$ is taken to be the unit outward vector normal to $\partial \Omega$. Let $e$ be an edge shared by two elements $T_m$ and $T_n$ of $\mathcal{T}_h$; we associate with $e$, once and for all, a unit normal vector $\bn_e$ directed from $T_m$ to $T_n$. We define the jump $[\psi]$ and the average $\{\psi\}$ of a function $\psi$ on $e$ by
	\[ [\psi]=(\psi|_{T_m})|_e-(\psi|_{T_n})|_e, \quad \{\psi\}=\frac{1}{2}(\psi|_{T_m})|_e+\frac{1}{2}(\psi|_{T_n})|_e.\] 
	If $e$ is adjacent to $\partial\Omega$, then  the jump and the average of $\psi$ on $e$ coincide with the value of $\psi$ on $e$.\\
	We define the following "broken" norm for  $l>0$ as
	\[\vertiii{\cdot}_l=\left(\sum_{T\in \mathcal{T}_h} \|\cdot\|^2_{l,T}\right)^{1/2}.\]
	We also need the following discontinuous spaces for  our subsequent analysis:
	\begin{align*}
		\bX &=\{\bw\in (L^2(\Omega))^2:\, \bw|_{T}\in \textbf{W}^{2,4/3}(T),\quad \forall T\in \mathcal{T}_h\},\\
		M &= \{q\in L^2(\Omega)/ \R:\, q|_T\in W^{1,4/3}(T),\quad\forall	T\in \mathcal{T}_h\}.
	\end{align*}
	We associate with the discontinuous spaces $\bX$ and $M$ the following norms:
	\begin{align*}
		\|\bv\|_\varepsilon &=(\vertiii{\nabla \bv}_0^2+J(\bv,\bv))^{1/2}\quad \forall \bv\in \bX,\\
		\|q\|_{M} &=\|q\|_{L^2(\Omega)/ \R}\quad \forall q\in M,
	\end{align*}
	where the jump term $J_0$ is defined as
	\[J_0(\bv,\bw)=\sum_{e\in\Gamma_h}\frac{\sigma_e}{|e|}\int_e[\bv]\cdot [\bw]\,ds.\]
	Here $|e|$ denotes the measure of the edge $e$ and $\sigma_e$ is a positive constant, defined for each edge $e$ and will be called penalty parameter.\\
	We recall the $L^p$-estimate for functions in $\bX$ in terms of the norm $\|\cdot\|_\varepsilon$ (\cite{GRW005}):  For each real number $p\in [2,\infty)$, there exists a constant $C(p)>0$ such that 
	\begin{equation}
		\|\bv\|_{L^p(\Omega)}\leq C(p)\|\bv\|_\varepsilon,\quad  \forall \bv\in \bX.\label{Lp}
	\end{equation} 
	We further introduce the bilinear forms $a:\textbf{X}\times\textbf{X}\rightarrow \mathbb{R}$ and $b:\textbf{X}\times M\rightarrow \mathbb{R}$ corresponding to the discontinuous Galerkin formulation as follows:
	\begin{align}
		a(\bw,\bv) &= \sum_{T\in \mathcal{T}_h}\int_T \nabla \bw:\nabla \bv\,dT -\sum_{e\in \Gamma_h} \int_e\{\nabla \bw\}\bn_e\cdot [\bv]\,ds +\epsilon\sum_{e\in \Gamma_h}\int_e \{\nabla \bv\}\bn_e\cdot [\bw]\,ds,\label{1.6}\\
		b(\bv,q) &= -\sum_{T\in \mathcal{T}_h}\int_T q\nabla \cdot \bv\,dT +\sum_{e\in \Gamma_h }\int_e \{q\}[\bv]\cdot \bn_e \,ds,\label{1.7}
	\end{align}
	where $\epsilon$ takes the constant value $1$ or $-1$.   The SIPG method is the case when  $ \epsilon=-1$ whereas the NIPG method is the case when  $ \epsilon=1$. We make the following assumptions throughout the paper:\\
	  If $\epsilon=-1$, then the penalty parameter $\sigma_e$ cannot be arbitrary. It must be bounded below by $\sigma_0>0$ and $\sigma_0$ is sufficiently large. If $\epsilon=1$, the  penalty parameter $\sigma_e$ can be simply equated to $1$ on all edges.
	\\ 
	We also define a trilinear form $c(\cdot,\cdot,\cdot)$ for the nonlinear convective term present in the system (\ref{8.1})-(\ref{8.4}) which is motivated from the Lesaint-Raviart upwinding scheme  (see \cite{LR74}), introduced in \cite{GRW005}.
	\begin{align}\label{4.8}
		c^\bw(\bu,\bz,\bt)=\sum_{T\in\mathcal{T}_h} & \left(\int_T(\bu\cdot \nabla \bz)\cdot \bt\,dT+{\int_{\partial T_-}|\{\bu\}\cdot \bn_T|(\bz^{int}-\bz^{ext})\cdot \bt^{int}}\,ds\right)\nonumber\\
		& +\frac{1}{2}\sum_{T\in \mathcal{T}_h}\int_T(\nabla\cdot \bu)\bz\cdot \bt\,dT -\frac{1}{2}\sum_{e\in \Gamma_h}\int_e[\bu]\cdot \bn_e\{\bz\cdot \bt\}\,ds ,\qquad \forall \bu,\bz,\bt\in \textbf{X},
	\end{align} 
	where 
	\[\partial T_-=\{\textbf{x}\in \partial T:\{\bw\}\cdot \bn_T<0\}.\]
	The superscript $\bw$ denotes the dependence of $\partial T_-$ on $\bw$ and the superscript int ( respectively ext) refers to the trace of the function on a side of $T$ coming from the interior  (respectively exterior) of $T$ on that side. When the side of $T$ belongs to $\partial \Omega$, then we take the exterior trace to be zero. The first two terms in the definition of $c(\cdot,\cdot,\cdot)$ were introduced in \cite{LR74} for solving transport problems; the last term is chosen so that $c$ satisfies (\ref{5.5}), which ensures its positivity. Note that, $c^\bw(\bu,\bz,\bt)$ can also be written as 
	\[c^\bw(\bu,\bz,\bt)=\sum_{T\in\mathcal{T}_h}  \left(\int_T(\bu\cdot \nabla \bz)\cdot \bt\,dT+{\int_{\partial T_-}|\{\bu\}\cdot \bn_T|(\bz^{int}-\bz^{ext})\cdot \bt^{int}}\,ds\right)-\frac{1}{2}b(\bu,\bz\cdot \bt). \]
	It is easy to see that, when $\bu,\bz,\bt\in H_0^1(\Omega)^2$, the trilinear form $c$ reduces to
	\begin{equation}
		c(\bu;\bz,\bt)=\int_\Omega (\bu\cdot \nabla \bz)\cdot \bt\,dT+\frac{1}{2}\int_\Omega (\nabla \cdot \bu)\bz\cdot \bt\,dT.\label{4.9}
	\end{equation}
	The superscript $\bw$ is dropped in (\ref{4.9}) since the integral on $\partial T_-$ disappears. It is proven in \cite{GRW005} that  $c(\cdot,\cdot,\cdot)$  satisfies the following ``integration by parts'' for all $\bu,\bz,\bt \in\textbf{X}$:
	\begin{align}
		c^{\bu}(\bu,\bz,\bt)=& -\sum_{T\in \mathcal{T}_h}\left(\int_T (\bu \cdot \nabla \bt)\cdot \bz\,dT+\frac{1}{2}\int_T(\nabla\cdot \bu)\bz\cdot\bt\,dT\right)+\frac{1}{2}\sum_{e\in \Gamma_h}\int_e [\bu]\cdot \bn_e \{\bz\cdot\bt\}\,ds\nonumber\\
		&  -\sum_{T\in \mathcal{T}_h}\int_{\partial T_-}|\{\bu\}\cdot \bn_T|\bz^{ext}\cdot (\bt^{int}-\bt^{ext})\,ds +\int_{\Gamma_+}|\bu\cdot \bn|\bz\cdot \bt\,ds,\label{integration}
	\end{align}
	where $\Gamma_+$ is the subset of $\partial\Omega$ where $\bu\cdot \bn>0$. In particular, if $\bz=\bt$, then for $\bu,\bz \in\textbf{X}$,  we obtain 
	\begin{equation}
		c^\bu(\bu,\bz,\bz)\geq 0.\label{5.5}
	\end{equation}
	We now consider the DG formulation of (\ref{8.1})-(\ref{8.4}):  Find the pair $(\bu (t),p(t))\in \bX \times M,~ t>0$, such that 
	\begin{align}
		(\bu _t(t),\bphi)+\mu (a(\bu(t),\bphi)+J_0(\bu(t),\bphi))+c^\bu(\bu(t),\bu(t),\bphi)+b(\bphi,p(t))&= (\textbf{f}(t),\bphi)\quad \forall \bphi\in \bX,\label{8.5}\\
		b(\bu(t),q)&= 0\quad \forall q\in M,\label{8.6}\\
		(\bu(0),\bphi)&= (\bu_0,\bphi)\quad \forall \bphi\in \bX.\label{8.7}
	\end{align}
	For the consistency of the scheme (\ref{8.5})-(\ref{8.7}), one can refer to \cite{KR05} (Lemma 3.2). \\
	To carry out our analysis, we first recall some standard trace and inverse inequalities, which hold true on each element $T$ in $\mathcal{T}_h$, with diameter $h_T$ ( for a proof, please refer to \cite{GR79}):
	\begin{lemma}\label{trace}
		For every $T$ in $\mathcal{T}_h$, the following inequalities hold
		\begin{align}
			\|\bv\|_{0,e} &  \leq C(h_T^{-1/2}\|\bv\|_{0,T}+h_T^{1/2}\|\nabla\bv\|_{0,T})\quad \forall e\in \partial T,\,\, \forall\bv\in \bX,\label{T1}\\
			\|\nabla\bv\|_{0,e} & \leq C(h_T^{-1/2}\|\nabla\bv\|_{0,T}+h_T^{1/2}\|\nabla^2 \bv\|_{0,T})\quad \forall e\in \partial T,\,\,\forall \bv\in \bX,\label{T2}\\
			\|\bv\|_{L^4(e)} & \leq Ch_T^{-3/4}(\|\bv\|_{0,T}+h_T\|\nabla \bv\|_{0,T})\quad \forall e\in \partial T,\,\, \bv\in \bX.\label{T3}
		\end{align}
	\end{lemma}
	\noindent
	Further, we state the  regularity estimates which will be used in the subsequent error analysis \cite{HR82}.
	\begin{lemma} \label{exactpriori}
		Let the assumptions ({\bf A1}) and ({\bf A2}) hold. Then there exists a positive constant $C$ such that the weak solutions of (\ref{w1})-(\ref{w3}) satisfy the following regularity estimates, for $t>0$
		\begin{align}
			\sup_{0<t<\infty}\{\|\bu\|_2+\|\bu_t\|+\|p\|_{H^1/R}\} &\leq C, \label{priori1} \\
			\sigma^{-1}(t)\int_0^t e^{2\al s}(\|\bu\|_{2}^2+\|p\|_{1}^2)\,ds &\leq C, \label{priori2}\\
			\sup_{0<t<\infty}\tau(t)\|\bu_t\|_1^2+e^{-2\al t}\int_0^t \sigma(s)(\|\bu_t\|_{2}^2+\|p_t\|_{1}^2)\,ds &\leq C ,\label{priori3}
		\end{align}
		where $\tau(t)=min\{t,1\}, \quad \sigma(t)=\tau(t)e^{2\alpha t}$.
	\end{lemma}
	
	\section{Semidiscrete discontinuous Galerkin formulation}\label{s3}
	\se
	
	On the triangulation, defined in {\bf Section} \ref{s2},  let $\bX_h\subset \bX$ and $M_h\subset M$ be two finite-dimensional subspaces for approximating velocity and pressure, respectively.
	For any positive integer $k\ge 1$, we define them as follows:
	\begin{align*}
		\bX_h & = \{\textbf{v}\in L^2(\Omega)^2\, :\,\, \forall T \in \mathcal{T}_h, \,\, \textbf{v}\in (\mathbb{P}_k(T))^2\},\\
		M_h & = \{q\in  L^2(\Omega)/ \R\, :\, \forall T\in \mathcal{T}_h,\,\, q\in \mathbb{P}_{k-1}(T)\}
	\end{align*}
	We assume the following approximation properties for the spaces $\bX_h$ and $M_h$. We can construct an operator $r_h\in \mathcal{L}(L^2(\Omega)/ \R;M_h)$ (see \cite{GRW005}), such that for any $T\in \mathcal{T}_h$,
	\begin{equation}
		\forall q\in L^2(\Omega)/ \R,~~\forall z_h\in \mathbb{P}_{k-1}(T),\quad  \int_T  z_h(r_h(q)-q)\,dT=0,\label{2}
	\end{equation}
	and for any real number $s\in [0, k]$,
	\begin{equation}
		\forall q\in H^s(\Omega)\cap L^2(\Omega)/ \R,\quad \|q-r_h(q)\|_{0,T}\leq C h_T^s|q|_{s,T}. \label{2.1}
	\end{equation}
	And $\bX_h$ allows a projection operator $\textbf{R}_h$.
	\begin{lemma}\label{Rhp}
		For $\bX_h$, there exists an operator $\textbf{R}_h\in \mathcal{L}(H^1(\Omega)^2;\textbf{X}_h(\Omega))$, such that for any $T\in \mathcal{T}_h$,
		\begin{align}
			\forall \bv\in \bH^1(\Omega),\quad  \forall q_h\in \mathbb{P}_{k-1}(T),\quad &\int_T  q_h \nabla\cdot (\textbf{R}_h(\bv)-\bv)\,dT = 0, \label{2.2}\\
			\forall \bv\in \bH^1(\Omega),~ \forall e\in \Gamma_h,\quad \forall q_h\in \mathbb{P}_{k-1}(e)^2,\quad  &\int_e  q_h\cdot [\textbf{R}_h(\bv)]\,ds = 0,\label{2.3}\\
			\forall \bv\in \bH_0^1(\Omega),~ \forall e\in \partial \Omega,\quad  \forall q_h \in \mathbb{P}_{k-1}(e)^2,\quad &\int_e  q_h \cdot \textbf{R}_h(\bv)\,ds = 0,\label{2.4}\\
			\forall \bv\in  \bH^{k+1}(\Omega)\cap\bH_0^1(\Omega),~~ &\|\bR_h(\bv)-\bv\|_\varepsilon\leq C h^k|\bv|_{k+1,\Omega},\label{2.5}\\
			\forall s \in [1, k+1],~ \forall \bv\in \bH^s(\Omega),~~ &\|\bv-\textbf{R}_h(\bv)\|_{L^2(T)}\leq C h_T^s|\bv|_{s,\Delta_T}, \label{2.6}
		\end{align}
		where $\Delta_T$ is a suitable macro-element containing $T$.
	\end{lemma}
	\noindent  For $k=1,2$ and $3$, the existence of this operator and (\ref{2.2})-(\ref{2.4}) follows from \cite{CR73,FS83,CF89}. The bounds (\ref{2.5}) and (\ref{2.6}) are proved in \cite{GRW005} and \cite{GS03}, respectively.
	Recall that, the operator $\bR_h$ satisfies (see \cite{GRW005})
	\begin{equation}
		\forall v\in \bH_0^1(\Omega)^2, \quad \forall q \in M_h,\quad b(\textbf{R}_h(v)-v,q)=0.\label{3.9}
	\end{equation}
	Furthermore, $\bR_h$ satisfies the following stability property (see \cite{GRW05}): there exists a constant $C$, independent of $h$, such that 
	\begin{equation}
		\|\bR_h \bv\|_\varepsilon\leq C |\bv|_{H^1(\Omega)}, \quad \forall \bv\in \bH_0^1(\Omega).\label{Rh}
	\end{equation}
	We now define  the semi-discrete discontinuous Galerkin approximations for the equations (\ref{w1})-(\ref{w3}).
	For all $t>0$, we seek a discontinuous approximation $(\bu_h (t),p_h(t))\in \bX_h \times M_h$ such that 
	\begin{align}
		(\bu_{ht}(t),\bphi_h)+\mu (a(\bu_h(t),&\bphi_h)+J_0(\bu_h(t),\bphi_h))
		+c^{\bu_h}(\bu_h(t),\bu_h(t),\bphi_h)+b(\bphi_h,p_h(t))= (\textbf{f}(t),\bphi_h), \label{8.8}\\
		& b(\bu_h(t),q_h)= 0, \quad \mbox{and} \quad
		(\bu_h(0),\bphi_h)= (\bu_0,\bphi_h).\label{8.9}
	\end{align}
	for $(\bphi_h,q_h)\in (\bX_h,M_h)$. In order to consider a discrete space analogous to $ \bJ_1$, we define the space $\bV_h$ by:
	\[\textbf{V}_h=\{\bv_h\in \textbf{X}_h:\forall q_h\in M_h,\quad b(\bv_h,q_h)=0\}.
	\]
	Now, an equivalent formulation of (\ref{8.8})--(\ref{8.9}) on $\bV_h$  reads as follows:
	Find $\bu_h(t)\in \bV_h$, such that for $t>0$
	\begin{equation}
		(\bu_{ht},\bphi_h)+\mu(a(\bu_h,\bphi_h)+J_0(\bu_h,\bphi_h))+c^{\bu_h}(\bu_h,\bu_h,\bphi_h)=(\bf f, \bphi_h),~\forall\quad \bphi_h\in \bV_h. \label{u_h}
	\end{equation}
	Below in Lemma \ref{coer}, we state the ellipticity of bilinear form  $(a+J_0)(\cdot,\cdot)$, which can be used to prove the well-posedness of the above discrete system(s).
	\begin{lemma}[\cite{W78}]\label{coer}
		 Assume that $\sigma_e$ is sufficiently large. Then for $\epsilon=-1$, there is a constant $K>0$ independent of $h$ such that 
		\begin{equation*}
			\forall \bv_h\in \textbf{X}_h,\quad a(\bv_h,\bv_h)+J_0(\bv_h,\bv_h)\geq K\|\bv_h\|_\varepsilon ^2.
		\end{equation*} 
	\end{lemma}
	\noindent
	 \begin{Remark}\label{Remark0}
			Using the definition of $\|\cdot\|_\varepsilon$-norm, the above lemma can be proved easily for the nonsymmetric case ($\epsilon=1$) with $K=1$. 
	\end{Remark} 
	\noindent
	Moreover, for the boundedness of bilinear form $a(\cdot,\cdot)$, we have the following result \cite{R08}.
	\begin{lemma}\label{cont}
		There exists a constant $C>0$, independent of $h$, such that, for all $\bv_h,\bw_h\in \bX_h$,
		\begin{equation*}
			|a(\bv_h,\bw_h)|\leq C \|\bv_h\|_\varepsilon\|\bw_h\|_\varepsilon.
		\end{equation*}
	\end{lemma}
	\noindent
	In Lemma \ref{inf-sup} below, we state a uniform discrete inf-sup condition for the pair of discontinuous spaces ($\widetilde{\bX}_h,M_h$), where
	\[ \widetilde{\textbf{X}}_h=\{\bv_h\in \textbf{X}_h: \forall e \in \Gamma_h, \quad \int_e \bq_h\cdot [\bv_h]\,ds=0,\quad \forall \bq_h\in (\mathbb{P}_{k-1}(e))^2\}. \]
	\begin{lemma}[\cite{GRW005}]\label{inf-sup}
	 There exist a constant $\beta^*>0$, independent of $h$, such that 
		\begin{equation*}
			\inf_{p_h\in M_h}\sup_{\bv_h\in \widetilde{\textbf{X}}_h}\frac{b(\bv_h,p_h)}{\|\bv_h\|_\varepsilon \|p_h\|_0}\geq \beta ^*.
		\end{equation*}
	\end{lemma}
	\noindent
	\noindent
	Now from the coercivity result  in Lemma \ref{coer}, the positivity (\ref{5.5}) and the inf-sup condition in  Lemma \ref{inf-sup}, the existence and uniqueness of the discrete Navier-Stokes solution to (\ref{8.8})-(\ref{8.9}) easily follow from \cite{KR05}. \\
	\noindent
	Before we proceed with our analysis, we recall trace inequalities for discrete space $\bX_h$, that are analogous to those from Lemma \ref{trace}.
	\begin{lemma}\label{distrace}
		For every element $T$ in $\mathcal{T}_h$, the following inequalities hold
		\begin{align}
			\|\bv_h\|_{0,e}&\leq Ch_T^{-1/2}\|\bv_h\|_{0,T}\quad \forall e\in \partial T, \,\,\forall \bv_h\in\bX_h,\label{T4}\\
			\|\nabla\bv_h\|_{0,e}&\leq Ch_T^{-1/2}\|\nabla\bv_h\|_{0,T}\quad \forall e\in \partial T, \,\,\forall \bv_h\in\bX_h,\label{T5}\\
			\|\nabla\bv_h\|_{0,T}&\leq Ch_T^{-1}\|\bv_h\|_{0,T}\quad \forall \bv_h\in\bX_h,\label{T6}\\
			\|\bv_h\|_{L^4(T)}&\leq Ch_T^{-1/2}\|\bv_h\|_{0,T}\quad \forall \bv_h\in\bX_h,\label{T7}
		\end{align}
		where $C$ is a constant independent of $h_T$ and $\bv_h$. 
	\end{lemma}
 Now, in Theorem \ref{semierror}, we state one of the main results of this article, which is related to the semi-discrete velocity error estimates.
	\begin{theorem}\label{semierror}
	Suppose the assumptions ({\bf A1}) and ({\bf A2}) hold. Further, let the discrete initial velocity $\bu_{h}(0)\in \bV_h$ with $\bu_{h}(0) =
		\bP_h\bu_0$, where $u_0\in {\bf H}^2\cap \bJ_1$. Then, there exists a positive constant $C$, independent of $h$, such that
		\begin{align*}
			\|(\bu-\bu_h)(t)\|+ h\|(\bu-&\bu_h)(t)\|_\varepsilon\leq Ce^{Ct}h^2.
		\end{align*}
	\end{theorem}
	\noindent
 Next, we aim to the derivation of results, which will lead to the proof of Theorem \ref{semierror}.\\
 
	\noindent
 {\bf\large{Approximation operators}} \vspace*{0.5cm}
		
		As mentioned in the introduction, we feel the need for appropriate approximation operators on the broken Sobolev spaces, which would allow us to obtain an optimal $L^{\infty}(\bL^2)$-norm error estimate for the discrete velocity, which are missing from the DG literature. Since we carry out our analysis for weakly divergence-free spaces, below, we derive new approximation properties for the space $ \bJ_1$.
	\begin{lemma}\label{ih}
		For every $\bv\in  \bJ_1\cap \bH^2$, there exists an approximation $i_h \bv\in \bV_h$ such that
		\begin{equation*}
			\|\bv-i_h \bv\|_\varepsilon\leq Ch\|\bv\|_2.
		\end{equation*}
		\begin{proof}
			From Lemma \ref{Rhp}, we have operator $\bR_h:\bH_0^1\rightarrow \bX_h$ satisfying $\|\bR_h \bv\|_\varepsilon\leq C|\bv|_1$ (see (\ref{Rh})) and 
			\[b(\bv-\bR_h \bv,q_h)=0, \quad \forall ~~ \bv\in \bH_0^1(\Omega), ~ q_h\in M_h.\]
			Restricting $\bR_h$ to $ \bJ_1\cap \bH^2$, we observe that
			$b(\bR_h \bv,q_h)=0$,  meaning
			$\bR_h\bv\in \bV_h$. We define this restriction as $i_h: \bJ_1\cap \bH^2\rightarrow \bV_h$. And from (\ref{2.5}), we have 
			\[\|\bv-i_h\bv\|_\varepsilon=\|\bv-\bR_h\bv\|\varepsilon\leq Ch\|\bv\|_2.\]
		\end{proof}
	\end{lemma}
	\noindent We now define a projection $\bP_h:\bL^2\rightarrow \bV_h$ which satisfies for $\bv\in\bL^2$
	\[(\bv-\bP_h\bv,\bs_h)=0,\quad \forall \bs_h\in \bV_h. \]
	The following lemma
	is a standard consequence of Lemma \ref{ih}.
	\begin{lemma}\label{P_h}
		The $L^2$- projection $\bP_h$ satisfies the following estimates:
		\begin{equation*}
			\|\bv-\bP_h\bv\|+h\|\bv-\bP_h\bv\|_\varepsilon\leq Ch^2|\bv|_2,\quad \forall\,\, \bv\in   \bJ_1\cap\bH^2.
		\end{equation*}
	\end{lemma}
	\noindent
	Finally, we define an approximation operator, a modified Stokes projection $\bS_h\bu\in \bV_h$, for the weak solution $\bu$ of the problem (\ref{w1})-(\ref{w3}), satisfying,
	\begin{align}
		\mu(a(\bu-\bS_h\bu,\bphi_h)+J_0(\bu-\bS_h\bu,\bphi_h))=-b(\bphi_h,p),\quad \forall \bphi_h\in \bV_h.\label{9.3}
	\end{align}
	\begin{lemma}\label{zetaerr}
		Let the asssumption $({\bf A2})$ hold true. Then $\bu-\bS_h\bu$ satisfies the following estimates:
		\begin{align}
			\|\bu-\bS_h\bu(t)\|^2+h^2\|\bu-\bS_h\bu(t)\|^2_\varepsilon\leq &  Ch^{4}\left( |\bu|_{2}^2+\frac{1}{\mu^2}|p|_1^2\right),\label{9.4}\\
			\|(\bu-\bS_h\bu)_t(t)\|^2+h^2\|(\bu-\bS_h\bu)_t(t)\|^2_\varepsilon \leq &  Ch^{4}\left( |\bu_t|_{2}^2+\frac{1}{\mu^2}|p_t|_1^2\right),\label{9.5}
		\end{align}
		where $C$ is a positive constant independent of $h$.
	\end{lemma}
	\begin{proof}
		Since
		\begin{align}\label{bzt01}
			\|\bu-\bS_h\bu\|_\varepsilon\leq \|\bu-\bP_h \bu\|_\varepsilon+\|\bP_h \bu-\bS_h \bu\|_\varepsilon,
		\end{align}
		it is sufficient to estimate $\bP_h \bu-\bS_h \bu$. In order to do that we choose $\bphi_h=\bP_h(\bu-\bS_h\bu)=\bu-\bS_h\bu-(\bu-\bP_h \bu)$ in (\ref{9.3}) to observe that
		\begin{align}
			\mu(&a( \bP_h\bu  -\bS_h\bu, \bP_h \bu-\bS_h\bu)+J_0(\bP_h \bu-\bS_h\bu,\bP_h \bu-\bS_h\bu))\nonumber\\
			&=-\mu a(\bu-\bP_h\bu,\bP_h \bu-\bS_h\bu)-\mu J_0(\bu-\bP_h \bu,\bP_h \bu-\bS_h\bu)-b(\bP_h\bu-\bS_h\bu,p-r_h(p)).\label{zetaener}
		\end{align}
		We expand the first term of right hand side to write as
		\begin{align}
			\mu a(\bu-\bP_h\bu,\bP_h\bu-\bS_h\bu)= &\mu\sum_{T\in \mathcal{T}_h}\int_T \nabla (\bu-\bP_h\bu):\nabla (\bP_h\bu-\bS_h\bu)\,dT\nonumber \\
			-\mu\sum_{e\in \Gamma_h} \int_e\{\nabla (\bu-\bP_h\bu)\}\bn_e\cdot & [\bP_h\bu-\bS_h\bu]\,ds + \mu\epsilon\sum_{e\in \Gamma_h}\int_e \{\nabla (\bP_h\bu-\bS_h\bu)\}\bn_e\cdot [\bu-\bP_h\bu]\,ds\nonumber\\
			=& S_1+S_2+S_3.\label{z1}
		\end{align}
		Using the Cauchy-Schwarz inequality, Young's inequality and Lemma \ref{P_h}, we obtain
		\begin{align*}
			|S_1|\leq \mu\sum_{T \in\mathcal{T}_h}\|\nabla(\bu-\bP_h\bu)\|_{L^2(T)}\|\nabla(\bP_h\bu-\bS_h\bu)\|_{L^2(T)}\leq \frac{K\mu}{10}\|\bP_h\bu-\bS_h\bu\|_\varepsilon^2+Ch^{2}|\bu|_{2}^2.
		\end{align*}
		For $S_2$, if the edge $e$ belongs to the element $T$, then by using trace inequality (\ref{T2}), we have 
		\[\bigg|\int_e\{\nabla (\bu-\bP_h\bu)\}\bn_e\cdot [\bP_h\bu-\bS_h\bu]\,ds\bigg|\leq C\Big(\|\nabla(\bu-\bP_h\bu)\|_{L^2(T)}+h_T\|\nabla^2(\bu-\bP_h\bu)\|_{L^2(T)}\Big)\frac{1}{|e|^{1/2}}\| [\bP_h\bu-\bS_h\bu]\|_{L^2(e)}. \]
		Let $L_h(\bu)$ denote the standard Lagrange interpolant of degree $1$. Then, by using inverse inequality (\ref{T6}), we obtain
		\begin{align*}
			\|\nabla^2(\bu-\bP_h\bu)\|_{L^2(T)}&\leq \|\nabla^2(\bu-L_h(\bu))\|_{L^2(T)}+\|\nabla^2(L_h(\bu)-\bP_h\bu)\|_{L^2(T)}\\
			&\leq \|\nabla^2(\bu-L_h(\bu))\|_{L^2(T)}+Ch_T^{-1}\|\nabla(L_h(\bu)-\bP_h\bu)\|_{L^2(T)}.
		\end{align*}
		Now Lemma \ref{P_h}, the standard approximation properties of $L_h$, the triangle inequality and the Cauchy-Schwarz inequality yield
		\[|S_2|\leq C\mu hJ_0(\bP_h\bu-\bS_h\bu,\bP_h\bu-\bS_h\bu)^{1/2}|\bu|_{2}\leq  \frac{K\mu}{10}\|\bP_h\bu-\bS_h\bu\|_\varepsilon^2+Ch^{2}|\bu|_{2}^2.\]
		Furthermore, trace inequality (\ref{T5}) and Lemma \ref{P_h} yield
		\begin{align*}
			|S_3|\leq C\mu\bigg(\sum_{e\in\Gamma_h}\frac{|e|}{\sigma_e}\|\{\nabla(\bP_h\bu-\bS_h\bu)\}\|_{L^2(e)}^2\bigg)^{1/2}\bigg(\sum_{e\in\Gamma_h}\frac{\sigma_e}{|e|}\|[\bu-\bP_h\bu]\|_{L^2(e)}^2\bigg)^{1/2}\leq \frac{K\mu}{10}\|\bP_h\bu-\bS_h\bu\|_\varepsilon^2+Ch^{2}|\bu|_{2}^2. 
		\end{align*}
		Using Cauchy-Schwarz's inequality, the jump term is bounded by virtue of Lemma \ref{P_h} as  follows:
		\begin{align}
			\mu|J_0(\bu-\bP_h\bu,\bP_h\bu-\bS_h\bu)|&\leq \mu J_0(\bu-\bP_h\bu,\bu-\bP_h\bu)^{1/2}J_0(\bP_h\bu-\bS_h\bu,\bP_h\bu-\bS_h\bu)^{1/2}\nonumber\\
			&\leq  \frac{K\mu}{10}\|\bP_h\bu-\bS_h\bu\|_\varepsilon^2+Ch^2|\bu|_{2}^2.\label{z2}
		\end{align}
		Owing to (\ref{2}), the pressure term is reduced to
		\[b(\bP_h\bu-\bS_h\bu,p-r_hp)=\sum_{e\in \Gamma_h}\int_e\{p-r_hp\}[\bP_h\bu-\bS_h\bu]\cdot \bn_e\,ds,\]
		which is bounded by using the Cauchy-Schwarz inequality, trace inequality (\ref{T1}) and the approximation result (\ref{2.1}) as  follows:
		\begin{align}\label{z3}
			|b(\bP_h\bu-\bS_h\bu,p-r_hp)| \leq & C\sum_{T\in \mathcal{T}_h}(\|p-r_hp\|_{0,T}+h_T\|\nabla(p-r_hp)\|_{0,T})J_0(\bP_h\bu-\bS_h\bu,\bP_h\bu-\bS_h\bu)^{1/2}\nonumber\\
			\leq & \frac{K\mu}{10}\|\bP_h\bu-\bS_h\bu\|_\varepsilon^2+\frac{Ch^{2}}{\mu}|p|_1^2.
		\end{align}
		By  incorporating Lemma \ref{coer}, (\ref{z1}), (\ref{z2}) and (\ref{z3}) in (\ref{zetaener}), we obtain 
		\begin{align}\label{st1}
			\mu K\|\bP_h\bu-\bS_h\bu \|_\varepsilon^2\leq   Ch^{2}|\bu|_{2}^2
			+\frac{Ch^{2}}{\mu}|p|_1^2.
		\end{align}
		Now  from (\ref{bzt01}), we complete the energy norm estimate of $\bu-\bS_h\bu$:
		\begin{equation}
			\|\bu-\bS_h\bu\|_\varepsilon^2\leq Ch^{2}\left( |\bu|_{2}^2+\frac{1}{\mu^2}|p|_1^2\right).\label{9.7}
		\end{equation}
		For $L^2$-norm estimate, we employ the Aubin-Nitsche duality argument. For fixed $h$, let $\{\bw,q\}$ be the pair of unique solution of the following steady state Stokes system:
		\begin{align}\label{9.9}
			-\mu \Delta \bw+\nabla q & =\bu-\bS_h\bu \quad\text{in $\Omega$},\quad\nabla\cdot \bw = 0\quad \text{in $\Omega$},\quad\bw|_{\partial \Omega} =0.
		\end{align}
		The above pair satisfies the following regularity result
		\begin{equation}
			\|\bw\|_2+\|q\|_1\leq C\|\bu-\bS_h\bu\|.\label{10.0}
		\end{equation}
		Now form $L^2$ inner product between (\ref{9.9}) and $\bu-\bS_h\bu$, and using the regularity of $\bw$ and $q$, we obtain
		\begin{align*}
			\|\bu-\bS_h\bu\|^2 = &\mu\sum_{T\in\mathcal{T}_h} \int_T\nabla \bw:\nabla(\bu-\bS_h\bu)\,dT- \mu\sum_{T\in\mathcal{T}_h}\int_{\partial T}(\nabla \bw \bn_T)\cdot (\bu-\bS_h\bu)\,ds \nonumber\\
			&-\sum_{T\in\mathcal{T}_h}\int_T q\nabla\cdot (\bu-\bS_h\bu)\,dT+\sum_{T\in\mathcal{T}_h}\int_{\partial T}q\bn_T\cdot (\bu-\bS_h\bu)\,ds\nonumber\\
			= & \mu \sum_{T\in\mathcal{T}_h} \int_T\nabla (\bu-\bS_h\bu):\nabla \bw\,dT - \mu\sum_{e\in \Gamma_h}\int_e \{\nabla \bw\}\bn_e\cdot [\bu-\bS_h\bu]\,ds+b(\bu-\bS_h\bu,q).
		\end{align*}
		We then use (\ref{9.3}) with $\bphi_h$ replaced by $\bP_h\bw$, and noting that $[\bw]\cdot \bn_e=0$ on each interior edge to obtain
		\begin{align}\label{10.1}
			\|\bu-\bS_h\bu(t)\|^2=&\mu \sum_{T\in\mathcal{T}_h} \int_T\nabla (\bu-\bS_h\bu):\nabla (\bw-\bP_h\bw)\,dT+\mu\epsilon\sum_{e\in \Gamma_h}\int_e \{\nabla (\bw-\bP_h\bw)\}\bn_e\cdot [\bu-\bS_h\bu]\,ds\nonumber \\
			&- \mu(1+\epsilon)\sum_{e\in \Gamma_h}\int_e \{\nabla \bw\}\bn_e\cdot [\bu-\bS_h\bu]\,ds+\mu\sum_{e\in \Gamma_h}\int_e \{\nabla (\bu-\bS_h\bu)\}\bn_e\cdot [\bP_h\bw-\bw]\,ds\nonumber\\
			&+\mu J_0(\bu-\bS_h\bu,\bw-\bP_h \bw)+b(\bu-\bS_h\bu,q)-b(\bP_h \bw-\bw,p-r_h(p)).
		\end{align} 
		Consider the SIPG form of $a(\cdot,\cdot)$ {\it i.e.} $\epsilon=-1$. Then the third term on the right hand side of (\ref{10.1}) will vanish.  Similar to the first and second  terms on the right hand side in (\ref{zetaener}), we bound the following terms and then using Lemma \ref{P_h}, (\ref{9.7}) and (\ref{10.0}) to find that
		\begin{align}
			\bigg|\mu \sum_{T\in\mathcal{T}_h} \int_T\nabla (\bu-\bS_h\bu):\nabla (\bw-\bP_h\bw)\,dT-\mu \sum_{e\in \Gamma_h}\int_e \{\nabla (\bw-\bP_h\bw)\}\bn_e\cdot [\bu-\bS_h\bu]\,ds\nonumber \\
			+\mu\sum_{e\in \Gamma_h}\int_e \{\nabla (\bu-\bS_h\bu)\}\bn_e\cdot [\bP_h\bw-\bw]\,ds +\mu J_0(\bu-\bS_h\bu,\bw-\bP_h \bw)\bigg|\nonumber\\
			\leq  Ch\|\bw\|_2\|\bu-\bS_h\bu\|_\varepsilon+ Ch^{2}|\bu|_{2}\|\bw\|_2\leq  \frac{1}{6}\|\bu-\bS_h\bu\|^2+Ch^{4}\big( |\bu|^2_{2}+\frac{1}{\mu^2}|p|^2_1\big).\label{z4}
		\end{align}
		And for the sixth term on the right-hand side of (\ref{10.1}) we have
		\begin{align}
			b(\bu-\bS_h\bu,q)&= b(\bu-\bS_h\bu-\bP_h\bu+\bS_h\bu,q)+b(\bP_h\bu-\bS_h\bu,q)= b(\bu-\bP_h \bu,q)+b(\bP_h\bu-\bS_h\bu,q-r_h(q)) \nonumber \\
			&= -\sum_{T\in \mathcal{T}_h}\int_T q\nabla\cdot (\bu-\bP_h\bu)\,dT+\sum_{e\in \Gamma_h}\int_e\{q\}[\bu-\bP_h \bu]\cdot \bn_e\,ds+b(\bP_h\bu-\bS_h\bu,q-r_h(q)) \label{z5}.
		\end{align}
		Furthermore, using integration by parts formula to the first term on the right hand side of (\ref{z5}) and noting that $q$ is continuous, we arrive at
		\begin{align*}
			b(\bu-\bS_h\bu,q)= \sum_{T\in\mathcal{T}_h}\int_T \nabla q\cdot (\bu-\bP_h\bu)\,dT+b(\bP_h\bu-\bS_h\bu,q-r_h(q)). 
		\end{align*}
		From Cauchy-Schwarz's inequality, Lemmas \ref{P_h} and \ref{trace}, (\ref{2.1}), (\ref{10.0}) and (\ref{st1}), we obtain 
		\begin{align}
			|b(\bu-\bS_h\bu,q)|  &\leq \bigg| Ch^{2}|q|_1|\bu|_{2}-\sum_{T\in\mathcal{T}_h}\int_T\nabla \cdot (\bP_h\bu-\bS_h\bu)(q-r_h(q))\,dT+\sum_{e\in\Gamma_h}\int_e\{q-r_h(q)\}[\bP_h\bu-\bS_h\bu]\cdot \bn_e\,ds\bigg|\nonumber\\
			& \leq Ch^{2}|\bu|_{2}\|\bu-\bS_h\bu\|+Ch|q|_1\|\bP_h\bu-\bS_h\bu\|_\varepsilon\leq \frac{1}{6}\|\bu-\bS_h\bu\|^2+Ch^{4}\left( |\bu|^2_{2}+\frac{1}{\mu^2}|p|^2_1\right).\label{z6}
		\end{align}
		Similarly, using Cauchy-Schwarz's inequality, we arrive at
		\begin{align}
			|b(\bP_h \bw-\bw,p-r_h(p))|
			\leq  Ch^{2}|p|_1\|\bw\|_2\leq\frac{1}{6}\|\bu-\bS_h\bu\|^2+Ch^{4}|p|_1^2 .\label{z7}
		\end{align}
		In view of (\ref{z4}), (\ref{z6}) and (\ref{z7}) in (\ref{10.1}), we complete the estimate (\ref{9.4}). \\
		Repeating the above set of arguments we arrive at the estimates (\ref{9.5}) involving $(\bu-\bS_h\bu)_t$. The only differences are instead of the equation (\ref{9.3}), we use the one obtained from differentiating in time, use $\bphi_h=\bP_h(\bu-\bS_h\bu)_t$ in it and finally for the dual problem, we take the right hand side as $(\bu-\bS_h\bu)_t$.
		This completes the proof of Lemma \ref{zetaerr}.
	\end{proof}
	\begin{Remark}\label{Remark1}
		In the case of NIPG formulation {\it i.e.} $\epsilon=1$, the third term on the right hand side of (\ref{10.1}) is nonzero. And here we will lose a power of $h$. Using Cauchy-Schwarz's inequality, Young's inequalty, trace inequality (\ref{T2}), and estimates (\ref{9.7}) and (\ref{10.0}), we can show that
		\begin{align*}
			\bigg|\sum_{e\in \Gamma_h}\int_e \{\nabla \bw\}\bn_e\cdot [\bu-\bS_h\bu]\,ds\bigg|&\leq C\bigg(\|\nabla \bw\|^2+\sum_{T\in\mathcal{T}_h}h_T^2\|\nabla^2\bw\|_{L^2(T)}^2\bigg)^{1/2}J_0(\bu-\bS_h\bu,\bu-\bS_h\bu)^{1/2}\\
			& \leq   C\|\bw\|_2\|\bu-\bS_h\bu\|_\varepsilon  \leq  \frac{1}{6}\|\bu-\bS_h\bu\|^2+Ch^{2}\big( |\bu|^2_{2}+\frac{1}{\mu^2}|p|^2_1\big).
		\end{align*}
		Thus, for the NIPG case the estimates (\ref{9.4}) and (\ref{9.5}) become
		\begin{align*}
			\|\bu-\bS_h\bu(t)\|^2+ \|\bu-\bS_h\bu(t)\|^2_\varepsilon\leq &  Ch^{2}\left( |\bu|_{2}^2+\frac{1}{\mu^2}|p|_1^2\right), \\
			\|(\bu-\bS_h\bu)_t(t)\|^2+\|(\bu-\bS_h\bu)_t(t)\|^2_\varepsilon \leq &  Ch^{2}\left( |\bu_t|_{2}^2+\frac{1}{\mu^2}|p_t|_1^2\right). 
		\end{align*}
	\end{Remark}
	\noindent
	Before we proceed to the next section,  we state some estimates of $\bu_h$. The estimates can be easily obtained using (\ref{8.8})-(\ref{8.9}).
	\begin{lemma}\label{uhpriori}
		Let the assumptions ({\bf A1}) and ({\bf A2}) hold. Then the semi-discrete discontinuous Galerkin approximation $\bu_h$ of the velocity $\bu$ satisfies, for $t>0$,
		\begin{align}
			\|\bu_h(t)\|+e^{-2\al t}\int_0^t e^{2\al s}\|\bu_h(s)\|_\varepsilon^2\,ds \leq C,\label{uh1}\\
			e^{-2\al t}\int_0^t e^{2\al s}\|\bu_{ht}(s)\|_\varepsilon^2\,ds \leq C,\label{uh3}\\
			e^{-2\al t}\int_0^t e^{2\al s}\|\bu_{htt}(s)\|_{-1,h}^2\,ds \leq C,\label{uh2}
		\end{align}
		where \[ \|\bu_{htt}\|_{-1,h}=\sup\bigg\{\frac{\langle\bu_{htt},\bphi_h\rangle}{\|\bphi_h\|_\varepsilon},~\bphi_h\in \bX_h,\bphi_h\neq 0\bigg\}. \]
		Moreover,
		\begin{align}
			\limsup_{t\rightarrow\infty}\|\bu_h(t)\|_\varepsilon\leq \frac{\|\brf\|_{L^\infty(\bL^2)}}{K\mu},\label{uh4}
		\end{align}
		where $C>0$ is a constant, independent of $h$, depends only on the given data.
	\end{lemma}
	\noindent
	 And estimates of the trilinear form $c(\cdot,\cdot,\cdot)$ which will be useful for our error analysis \cite{GRW05}.
	\begin{lemma}
		(i) Assume that $\bu\in \bW^{1,4}(\Omega)$. There exists a positive constant $C$ independent of $h$ such that 
		\begin{align}
			|c(\bv_h,\bu,\bw_h)|\leq C\|\bv_h\||\bu|_{\bW^{1,4}}\|\bw_h\|_\varepsilon, \quad \forall \bv_h,\bw_h\in\bV_h.\label{tri1}
		\end{align}
		(ii) For any $\bv \in \bX$, $\bv_h$, $\bw_h$ and $\bz_h$ in $\bX_h$, we have the following estimate:
		\begin{align}
			|c^{\bv}(\bv_h,\bw_h,\bz_h)|\leq C\|\bv_h\|_\varepsilon\|\bw_h\|_\varepsilon\|\bz_h\|_\varepsilon.\label{tri2}
		\end{align} 
	\end{lemma}
	\section{Error estimates for velocity}\label{s4}
	\se
	In this Section, we discuss optimal error estimates for the error $\de=\bu-\bu_h$. We split the error into two parts, $\de=\bxi+\be$, where $\bxi=\bu-\bv_h$ represents the error inherent in the DG finite element approximation of a linearized (Stokes) problem, and $\be=\bv_h-\bu_h$ represents the error caused by the presence of the nonlinearity in problem (\ref{8.1}). The linearized equation to be satisfied by the auxiliary function $\bv_h$ is:
	\begin{align}
		(\bv_{ht},\bphi_h)+\mu(a(\bv_h,\bphi_h)+J_0(\bv_h,\bphi_h)) = (\brf,\bphi_h)-c^\bu(\bu,\bu,\bphi_h)\quad \forall \bphi_h\in \bV_h.\label{9.1}
	\end{align}  
	Below we derive some estimates of $\bxi$. Subtracting (\ref{9.1}) from (\ref{8.5}), the equation in $\bxi$ can be written as 
	\begin{align}
		(\bxi_t,\bphi_h)+\mu(a(\bxi,\bphi_h)+J_0(\bxi,\bphi_h))=-b(\bphi_h,p), \quad \bphi_h\in \bV_h.\label{9.2}
	\end{align}
	\begin{lemma}\label{xiL2err}
		Suppose the assumptions ({\bf A1}) and ({\bf A2}) hold. Let $\bv_h(t) \in \bV_h$ be a solution of (\ref{9.1}) with initial condition $\bv_h(0)=\bP_h\bu_0$. Then, for $0\leq t<T$, there is a positive constant $C$, independent of $h$, such that 
		\begin{equation*}
			\int_0^te^{2\alpha s}\|\bxi(s)\|^2\, ds\leq Ch^4\sigma(t),
		\end{equation*}
		where $\tau(t)=min\{t,1\}, ~\sigma(t)=\tau(t)e^{2\alpha t}$ and  $0<\alpha\leq\frac{\mu K}{2C}$.
	\end{lemma}
	\begin{proof}
		Choosing $\bphi_h=\bP_h\bxi$ in (\ref{9.2}) and using Lemma \ref{coer}, we arrive at
		\begin{align}
			\frac{1}{2}\frac{d}{dt}\|\bxi\|^2+\mu K\|\bP_h\bxi\|^2_\varepsilon\leq (\bxi_t,\bu-\bP_h\bu)-\mu a(\bu-\bP_h\bu,\bP_h\bxi)-\mu J_0(\bu-\bP_h\bu,\bP_h\bxi)-b(\bP_h\bxi,p-r_hp).\label{13.1}
		\end{align}
		We first note that
		\[(\bxi_t,\bu-\bP_h\bu)=\frac{1}{2}\frac{d}{dt}\|\bu-\bP_h\bu\|^2.\] 
		Using Cauchy-Schwarz's inequality, the jump term is bounded by the virtue of Lemma \ref{P_h} as follows:
		\[\mu|J_0(\bu-\bP_h\bu,\bP_h\bxi)|\leq \mu J_0(\bu-\bP_h\bu,\bu-\bP_h\bu)^{1/2}J_0(\bP_h\bxi,\bP_h\bxi)^{1/2}\leq  \frac{K\mu}{6}\|\bP_h\bxi\|_\varepsilon^2+Ch^{2}|\bu|_{2}^2.\]
		The term $b(\bP_h\bxi,p-r_hp)$ can be bounded as the third term on the right hand side of (\ref{zetaener}) from Lemma {\ref{zetaerr}} as  follows:
		\begin{align*}
			|b(\bP_h\bxi,p-r_hp)|  \leq \frac{K\mu}{6}\|\bP_h\bxi\|_\varepsilon^2+Ch^{2}|p|_{1}^2.
		\end{align*}
		Finally, the second term on the right hand side in (\ref{13.1}) can be handled  as the first term on the right hand side of (\ref{zetaener}) in Lemma {\ref{zetaerr}} as follows:
		\begin{align*}
			\mu|a(\bu-\bP_h\bu,\bP_h\bxi)| \leq  \frac{K\mu}{6}\|\bP_h\bxi\|_\varepsilon^2+Ch^{2}|\bu|_{2}^2.
		\end{align*}
		Combining all the above estimates  in (\ref{13.1}) and using triangle inequality with Lemma \ref{P_h}, we obtain
		\begin{equation}
			\frac{d}{dt}\|\bxi\|^2+\mu K\|\bxi\|^2_\varepsilon\leq \frac{d}{dt}\|\bu-\bP_h\bu\|^2+Ch^{2}(|\bu|_{2}^2+|p|_1^2).\label{50.1}
		\end{equation}
		Multiplying (\ref{50.1}) by $e^{2\alpha t}$ and using the $L^p$-estimate (\ref{Lp}) with $p=2$, we find
		\begin{align*}
			\frac{d}{dt}(e^{2\alpha t}\|\bxi\|^2)+(\mu K-2\alpha C) e^{2\alpha t}\|\bxi\|^2_\varepsilon\leq \frac{d}{dt}(e^{2\alpha t}\|\bu-\bP_h\bu\|^2)+Ch^{2}e^{2\alpha t}(|\bu|_{2}^2+|p|_1^2).
		\end{align*}
		By setting $\alpha\leq\frac{\mu K}{2C}$, integrating from $0$ to $t$ and observing $\|\bxi(0)\|$ is of the order $h$, we obtain 
		\begin{align}
			\int_0^t e^{2\alpha s}\|\bxi(s)\|^2_\varepsilon\, ds\leq Ch^{2}\int_0^t e^{2\alpha s}(|\bu(s)|_{2}^2+|p(s)|_{1}^2)\,ds.\label{13.2}
		\end{align}

		\noindent
		To estimate $L^2$- norm error, we  use the following duality argument: For fixed $h > 0$  and $t\in (0,T)$, let $\bw(s)\in \bJ_1,\, q(s)\in H^1$, be the unique solution of the backward problem
		\begin{align}
			\bw_s+\mu\Delta \bw-\nabla q =e^{2\alpha s}\bxi, \quad 0\leq s\leq t,\label{10.7}
		\end{align}
		with $\bw(t)=0$ satisfying
		\begin{equation}
			\int_0^t e^{-2\alpha s}(\|\Delta\bw(s)\|^2+\|\bw_s(s)\|^2+\|\nabla q(s)\|^2)\, ds \leq C\int_0^t e^{2\alpha s}\|\bxi(s)\|^2\,ds.\label{back stokes bound}
		\end{equation}
		Form $L^2$-inner product between (\ref{10.7}) and $\bxi$ to obtain
		\begin{align}
			e^{2\alpha s}\|\bxi\|^2 &=(\bxi,\bw_s)-\displaystyle\sum_{T\in\mathcal{T}_h} \int_T\mu\nabla \bxi:\nabla \bw\,dT+ \sum_{T\in\mathcal{T}_h}\int_{\partial T}(\mu\nabla \bw \bn_T)\cdot \bxi\,ds +\sum_{T\in\mathcal{T}_h}\int_T q\nabla\cdot \bxi\,dT-\sum_{T\in\mathcal{T}_h}\int_{\partial T}q\bn_T\cdot \bxi\,ds,\nonumber\\
			& = (\bxi,\bw_s)-\displaystyle\sum_{T\in\mathcal{T}_h} \int_T\mu\nabla \bxi:\nabla \bw\,dT +\displaystyle\sum_{e\in \Gamma_h}\int_e \{\mu\nabla \bw\}\bn_e\cdot [\bxi]\,ds-b(\bxi,q).\label{10.9}
		\end{align}
		Using (\ref{9.2}) with $\bphi_h=\bP_h \bw$ and (\ref{10.9}), we obtain
		\begin{align}\label{13.5}
			e^{2\alpha s}\|\bxi\|^2= & (\bxi,\bw_s)+(\bxi_s,\bP_h \bw)-\sum_{T\in\mathcal{T}_h}\int_T \mu\nabla \bxi:\nabla(\bw-\bP_h \bw)\,dT +(1+\epsilon)\sum_{e\in\Gamma_h}\int_e\{\mu\nabla \bw\}\bn_e\cdot [\bxi]\,ds\nonumber\\
			& \quad -\epsilon \sum_{e\in \Gamma_h}\int_e\{\mu \nabla (\bw-\bP_h \bw)\}\bn_e\cdot [\bxi]\,ds-\sum_{e\in \Gamma_h}\int _e \{\mu \nabla \bxi\}\bn_e\cdot [\bP_h \bw]\,ds\nonumber\\
			& \quad +\mu J_0(\bxi,\bP_h \bw)-b(\bxi, q-r_h(q))+b(\bP_h \bw,p_h-p).
		\end{align}
		Consider $\epsilon=-1$. Using Cauchy-Schwarz's inequality, (\ref{T2}), Lemmas \ref{distrace} and \ref{P_h}, (\ref{2.1}) and the fact that $[\bw]=0$,  we  easily obtain
		\begin{align*}
			&\bigg|\sum_{T\in\mathcal{T}_h}\int_T \mu\nabla \bxi:\nabla(\bw-\bP_h \bw)\,dT  - \sum_{e\in \Gamma_h}\int_e\{\mu \nabla (\bw-\bP_h \bw)\}\bn_e\cdot [\bxi]\,ds\\
			& +\sum_{e\in \Gamma_h}\int _e \{\mu \nabla \bxi\}\bn_e\cdot [\bP_h \bw]\,ds -\mu J_0(\bxi,\bP_h \bw)+b(\bxi, q-r_h(q)) -b(\bP_h \bw,r_hp-p)\bigg|\\
			&\leq Ch(\|\bw\|_2+\|q\|_1)\|\bxi\|_\varepsilon+Ch^{2}|p|_1\|\bw\|_2.
		\end{align*}
		Using the definition of $\bP_h$, we rewrite
		\begin{align*}
			(\bxi,\bw_s)+(\bxi_s,\bP_h \bw)=\frac{d}{ds}(\bxi,\bw)-(\bxi_s,\bw-\bP_h \bw)=\frac{d}{ds}(\bxi,\bw)-\frac{d}{ds}(\bxi,\bw-\bP_h\bw)+(\bu-\bP_h\bu,\bw_s).
		\end{align*}
		Now from (\ref{13.5}), using Cauchy-Schwarz's and Young's inequalities, we  find
		\begin{align}
			e^{2\alpha s}\|\bxi\|^2\leq &\frac{d}{ds}(\bxi,\bP_h\bw)+ \delta e^{-2\alpha s}(\|\bw\|_2^2+\|\bw_s\|^2+\|q\|_1^2)+C\delta^{-1}h^2e^{2\alpha s}\|\bxi\|_\varepsilon\nonumber\\
			&+C\delta^{-1}h^{4}e^{2\alpha s}(|\bu|_{k+1}^2+|p|^2_k).\label{13.6}
		\end{align}
		On integrating (\ref{13.6}) with respect to $s$ from $0$ to $t$ and using (\ref{back stokes bound}), we obtain the following estimate
		\begin{align*}
			\int_0^t e^{2\alpha s}\|\bxi\|^2\,ds\leq & (\bxi(t),\bP_h\bw(t))- (\bxi(0),\bP_h\bw(0))+Ch^2\int_0^te^{2\alpha s}\|\bxi\|^2_\varepsilon\,ds+\delta\int_0^te^{2\alpha s}\|\bxi\|^2\, ds\\
			&+Ch^{4}\int_0^te^{2\alpha s}(|\bu|^2_{2}+|p|^2_1)\,ds.
		\end{align*}
		Choosing $\delta$ appropriately, we then apply the estimate (\ref{13.2}) and {\it a priori} estimate (\ref{priori2}) to complete the rest of the proof.
		\begin{Remark}\label{Remark2}
			For the NIPG case, similar to the Remark \ref{Remark1}, the fourth term on the right hand side  of (\ref{13.5}) can be bounded as
			\[\bigg|\sum_{e\in\Gamma_h}\int_e\{\mu\nabla \bw\}\bn_e\cdot [\bxi]\,ds\bigg|\leq   C\|\bw\|_2\|\bxi\|_\varepsilon\]
			which implies that
			\begin{equation*}
				\int_0^te^{2\alpha s}(\|\bxi(s)\|^2+\|\bxi(s)\|^2_\varepsilon)\, ds\leq Ch^2\sigma(t).
			\end{equation*}
			Thus in this case the estimate is subotimal.
		\end{Remark}
	\end{proof}
	\noindent
	For optimal error estimates of $\bxi$ in $L^\infty(\bL^2)$-norm, 
	we decompose it as follows:
	\[\bxi=(\bu-\bS_h\bu)+(\bS_h\bu-\bv_h)=\bzeta+\bt.\]
	Since the estimates of $\bzeta=\bu-\bS_h\bu$ are known from Lemma \ref{zetaerr}, it is sufficent to estimate $\bt$ which would allow us to draw the following conclusion.
	\begin{lemma} \label{xierror}
		There is a positive constant $C$, independent of $h$, such that for $t>0$, $\bxi$ satisfies the following estimates
		\begin{equation*}
			\|\bxi(t)\|+h\|\bxi(t)\|_\varepsilon\leq Ch^{2},\quad 0\leq t\leq T.
		\end{equation*}
	\end{lemma}
	\begin{proof}
		We consider the equation in $\bt$.
		\[(\bt_t,\bphi_h)+\mu(a(\bt,\bphi_h)+J_0(\bt,\bphi_h))=-(\bzeta_t,\bphi_h),\quad \forall \bphi_h\in\bV_h.\]
		Setting $\bphi_h=\bt$  in the above equation and using Lemma \ref{coer}, we  arrive at
		\begin{align*}
			\frac{1}{2}\frac{d}{dt}\|\bt\|^2+\mu K\|\bt\|^2_\varepsilon\leq -(\bzeta_t,\bt). 
		\end{align*}
		Multiply by $\sigma(t)$ and integrate the resulting  inequality with respect to time.
		\begin{align}
			\sigma(t)\|\bt\|^2+\mu K\int_0^t\sigma(s)\|\bt(s)\|^2_\varepsilon &\leq C\int_0^t (2\sigma_s(s)\|\bt(s)\|^2+\frac{\sigma^2(s)}{\sigma_s(s)}\|\bzeta_t\|^2)\,ds\nonumber\\
			& \leq C\int_0^te^{2\al s}(\|\bzeta(s)\|^2+\|\bxi(s)\|^2)\,ds+\int_0^t\frac{\sigma^2(s)}{\sigma_s(s)}\|\bzeta_t\|^2\,ds.\label{theta3}
		\end{align}
		Using the estimates (\ref{9.4}), (\ref{9.5}), Lemma \ref{xiL2err}, (\ref{priori2}) and (\ref{priori3}) in (\ref{theta3}), we  obtain
		\begin{align*}
			\|\bt(t)\|^2+\sigma^{-1}(t)\int_0^t\sigma(s)\|\bt(s)\|^2_\varepsilon &\leq Ch^4.
		\end{align*}
		By using the inverse relation (\ref{T6}), we  now conclude that
		\begin{align*}
			\|\bt(t)\|^2+h^2\|\bt(t)\|_\varepsilon^2 \leq Ch^4. 
		\end{align*}
		This along with (\ref{9.4}) and {\it a priori} estimate (\ref{priori2}) give us the desired result.
	\end{proof}
	 \begin{Remark}\label{Remark3}
	With the help of Remarks \ref{Remark0}, \ref{Remark1} and \ref{Remark2}, for NIPG case, we can show that
	\begin{equation*}
			\|\bxi(t)\|+h\|\bxi(t)\|_\varepsilon\leq Ch,\quad 0\leq t\leq T.
		\end{equation*}
	\end{Remark}
	\noindent  We are now left with the estimate of $\be=\bv_h-\bu_h$.
	\begin{lemma}\label{etaerror}
Suppose the assumptions ({\bf A1}) and ({\bf A2}) hold. Let $\bv_h(t) \in \bV_h$ be a solution of (\ref{9.1}) corresponding to the initial value $\bv_h(0)=\bP_h\bu_0$. Then, there exists a positive constant $C$, independent of $h$ such that for $0 \leq t <T$, the error $\be$, satisfies
		\begin{equation*}
			\| \be(t)\|+h\| \be(t)\|_\varepsilon\leq Ce^{Ct}h^2.
		\end{equation*}
	\end{lemma}
	\begin{proof}
		From the equations (\ref{9.1}) and (\ref{u_h}), satisfied by $\bv_h$ and $\bu_h$, respectively, we obtain 
		\[ (\be_t,\bphi_h)+\mu(a(\be,\bphi_h)+J_0(\be,\bphi_h))=c^{\bu_h}(\bu_h,\bu_h,\bphi_h)-c^\bu(\bu,\bu,\bphi_h) \quad \text{for $ \bphi_h\in \bV_h$}.\]
		Setting $\bphi_h=\be$ and using Lemma \ref{coer}, we find 
		\begin{equation}\label{10.6}
			\frac{1}{2}\frac{d}{dt}\|\be\|^2+\mu K\|\be\|_\varepsilon^2\leq c^{\bu_h}(\bu_h,\bu_h,\be)-c^\bu(\bu,\bu,\be). 
		\end{equation}
		We first note that since $\bu$ is continuous, we can rewrite
		\[c^\bu(\bu,\bu,\be)=c^{\bu_h}(\bu,\bu,\be).\]
		Secondly, whenever there is no confusion, we drop the superscript in the nonlinear terms.
		Let us  now rewrite the nonlinear terms  as follows:
		\begin{align}
			c(\bu_h,\bu_h,\be)-c(\bu,\bu,\be)& =-c(\bu,\bxi,\be)+c(\bxi,\bxi,\be)-c(\bxi,\bu,\be)+c(\be,\bxi,\be)-c(\be,\bu,\be)-c(\bu_h,\be,\be) \nonumber \\
			& \le -c(\bu,\bxi,\be)+c(\bxi,\bxi,\be)-c(\bxi,\bu,\be)+c(\be,\bxi,\be)-c(\be,\bu,\be). \label{c0}
		\end{align}
		The last term is non-negative and is dropped,  following (\ref{5.5}).  To bound the rest of the terms, we proceed as follows. A use of estimate (\ref{tri1}), Young's inequality and Sobolev's inequality implies
		\begin{align}
			|c(\be,\bu,\be)|\leq C \|\be\|_{L^2(\Omega)}|\bu|_{\bW^{1,4}(\Omega)}\|\be\|_\varepsilon\leq \frac{K\mu}{64}\|\be\|^2_{\varepsilon}+C\|\bu\|_2^2\|\be\|^2.\label{c4}
		\end{align}
		Using  Cauchy-Schwarz's inequality, Young's inequality, (\ref{Lp}) and Lemmas \ref{trace} and \ref{distrace}, the second and fourth nonlinear terms on the right hand side of (\ref{c0}) can be bounded as follows
		\begin{align}\label{c5}
			|c (\bxi,\bxi,\be)|= &\bigg|\sum_{T\in\mathcal{T}_h} \int_T (\bxi\cdot \nabla\bxi)\cdot \be\,dT+\sum_{T\in\mathcal{T}_h}\int_{\partial T_-}|\{\bxi\}\cdot \bn_T|(\bxi^{int}-\bxi^{ext})\cdot \be^{int}\,ds +\frac{1}{2}\sum_{T\in\mathcal{T}_h} \int_T(\nabla\cdot \bxi)\bxi\cdot \be\,dT \nonumber\\
			&-\frac{1}{2}\sum_{e\in\Gamma_h}\int_e[\bxi]\cdot \bn_e\{\bxi\cdot\be\}\,ds\bigg|\nonumber\\
			\leq  \sum_{T\in\mathcal{T}_h}\|\bxi\|_{L^4(T)}& \|\nabla\bxi\|_{L^2(T)}\|\be\|_{L^4(T)}+C\sum_{T\in\mathcal{T}_h}h_T^{-3/4}(\|\bxi\|_{L^2(T)}+h_T\|\nabla\bxi\|_{L^2(T)})|e|^{1/4}|e|^{-1/2}\|[\bxi]\|_{L^2(e)}\|\be\|_{L^4(T)}\nonumber\\
			+C\sum_{T\in\mathcal{T}_h}\|\nabla\bxi &\|_{L^2(T)}\|\bxi\|_{L^4(T)}\|\be\|_{L^4(T)}+C\sum_{T\in\mathcal{T}_h}|e|^{1/4}|e|^{-1/2}\|[\bxi]\|_{L^2(e)}h_T^{-3/4}(\|\bxi\|_{L^2(T)}+h_T\|\nabla\bxi\|_{L^2(T)})\|\be\|_{L^4(T)}\nonumber\\
			&\leq   \frac{K\mu}{64}\|\be\|^2_\varepsilon+C h^{-2}(\|\bxi\|^2+h^2\|\bxi\|_\varepsilon^2)\|\bxi\|_\varepsilon^2
		\end{align}
		and
		\begin{align}
			|c(\be,\bxi,\be)|= &\bigg|\sum_{T\in\mathcal{T}_h}\int_T (\be\cdot \nabla\bxi)\cdot \be\,dT+\sum_{T\in\mathcal{T}_h}\int_{\partial T_-}|\{\be\}\cdot \bn_T|(\bxi^{int}-\bxi^{ext})\cdot \be^{int}\,ds +\frac{1}{2}\sum_{T\in\mathcal{T}_h} \int_T(\nabla\cdot \be)\bxi\cdot \be\,dT \nonumber\\
			&-\frac{1}{2}\sum_{e\in\Gamma_h}\int_e[\be]\cdot \bn\{\bxi\cdot\be\}\,ds\bigg|\nonumber\\
			\leq  \sum_{T\in\mathcal{T}_h}\|\be\|_{L^4(T)}&\|\nabla \bxi\|_{L^2(T)}\|\be\|_{L^4(T)}+C \sum_{T\in\mathcal{T}_h}\|\be\|_{L^4(T)}\frac{1}{|e|^{1/2}}\|[\bxi]\|_{L^2(e)}\|\be\|_{L^4(T)}\nonumber\\
			+C\sum_{T\in\mathcal{T}_h}\|\nabla \be &\|_{L^2(T)}\|\bxi\|_{L^4(T)}\|\be\|_{L^4(T)} +C\sum_{T\in\mathcal{T}_h} h_T^{-3/4}\|\be\|_{L^2(T)} h_T^{-1/2}(\|\bxi\|_{L^2(T)}+h_T\|\nabla \bxi\|_{L^2(T)})h_T^{-1/4}\|\be\|_{L^4(T)}\nonumber\\
			& \leq \frac{K\mu}{64}\|\be\|^2_\varepsilon+C h^{-4}(\|\bxi\|^2+h^2\|\bxi\|_\varepsilon^2)\|\be\|^2.\label{c3}
		\end{align}
		For the first nonlinear term on the right hand side of (\ref{c0}), following (\ref{integration}), we rewrite it as
		\begin{align}
			c(\bu;\bxi,\be)=& -\sum_{T\in\mathcal{T}_h}\int_T(\bu\cdot \nabla \be)\cdot \bxi\,dT-\frac{1}{2}\sum_{T\in\mathcal{T}_h}\int_T(\nabla\cdot \bu)\bxi\cdot \be\,dT\nonumber\\
			& +\frac{1}{2}\sum_{e\in\Gamma_h}\int_e[\bu]\cdot \bn_e\{\bxi\cdot\be\}\,ds-\sum_{T\in \mathcal{T}_h}\int_{\partial T_-}|\bu\cdot \bn_T|\bxi^{ext}\cdot (\be^{int}-\be^{ext})\,ds\nonumber\\
			& +\int_{\Gamma_+}|\bu\cdot \bn|\bxi\cdot \be \,ds\nonumber\\
			&=A_1+A_2+A_3+A_4+A_5\label{c1}
		\end{align}
		Since $\bu$ is continuous, it is clear that $A_3=A_5=0$. For $A_1$ and $A_2$, we
		use the Cauchy-Schwarz's and  the Young's inequalities and the estimate (\ref{Lp}).
		\begin{align}
			|A_1|+|A_2|&\leq C \|\bu\|_{L^\infty(\Omega)}\|\nabla \be\|_{L^2(\Omega)}\|\bxi\|_{L^2(\Omega)}+ C\|\nabla \bu\|_{L^4(\Omega)}\|\bxi\|_{L^2(\Omega)}\|\be\|_{L^4(\Omega)} \nonumber\\
			&\leq C\|\bu\|_2\|\bxi\|\|\be\|_\varepsilon\leq C\|\bu\|_2^2\|\bxi\|^2+\frac{K\mu}{64}\|\be\|^2_\varepsilon \label{c11}
		\end{align}
		
		\noindent A use of (\ref{T1}) leads to the following bound of $A_4$:
		\begin{align}
			|A_4|& \leq C\|\bu\|_{L^\infty(\Omega)}\sum_{T\in \mathcal{T}_h}\|\bxi\|_{L^2(\partial T)}|e|^{1/2-1/2}\|[\be]\|_{L^2(\partial T)}\nonumber
			\\
			&\leq C\|\bu\|_2\sum_{T\in \mathcal{T}_h}|e|^{1/2}h_T^{-1/2}(\|\bxi\|_{L^2(T)}+h_T\|\nabla\bxi\|_{L^2(T)})\frac{1}{|e|^{1/2}}\|[\be]\|_{L^2(\partial T)}\nonumber\\
			& \leq C\|\bu\|_2^2(\|\bxi\|^2+h^2\|\bxi\|_\varepsilon^2)+\frac{K\mu}{64}\|\be\|^2_\varepsilon.\label{c13}
		\end{align}
		Finally, for the third nonlinear term on the right hand side of (\ref{c0}), we use $(\nabla\cdot \bu)v_i=\nabla\cdot (\bu v_i)-\bu\cdot \nabla v_i$ thereby
		\begin{align*}
			\int_T(\nabla \cdot \bu)\bv\cdot \bw\,dT = &\int_T \nabla\cdot (\bv\otimes \bu)\cdot \bw\,dT-\int_T(\bu\cdot \nabla \bv)\cdot \bw\,dT\\
			=&-\int_T(\bu\cdot \nabla \bw)\cdot \bv\,dT+\int_{\partial T}(\bu^{int} \cdot \bn_T)\bv^{int}\cdot \bw^{int}\,ds -\int_T(\bu\cdot \nabla \bv)\cdot \bw\,dT.
		\end{align*}
		This allows us the following reformulation.
		\begin{align}
			c(\bxi, \bu, \be)=& \frac{1}{2}\sum_{T\in\mathcal{T}_h}\int_T (\bxi\cdot \nabla\bu)\cdot \be\,dT- \frac{1}{2}\sum_{T\in\mathcal{T}_h}\int_T (\bxi\cdot \nabla \be)\cdot \bu\,dT\nonumber\\
			&- \sum_{T\in\mathcal{T}_h}\int_{\partial T_-} |\{\bxi\}\cdot \bn_T|(\bu^{int}-\bu^{ext})\cdot \be^{int}\,ds-\frac{1}{2}\sum_{e\in\Gamma_h}\int_e[\bxi]\cdot \bn_e\{\bu\cdot \be\}\,ds\nonumber\\
			&+ \frac{1}{2}\sum_{T\in \mathcal{T}_h}\int_{\partial T} (\bxi^{int}\cdot \bn_T)\bu^{int}\cdot \be^{int} \,ds\nonumber\\
			= & A_6+A_7+A_8+A_9+A_{10}.\label{c2}
		\end{align}
		Since $\bu$ is continuous, we have, $A_8=0$. The terms $A_{6}$ and $A_{7}$ are bounded using (\ref{Lp}) as follows:
		\begin{align}
			|A_{6}| & \leq \|\bxi\|_{L^2(\Omega)}\|\|\nabla\bu\|_{L^4(\Omega)}\|\be\|_{L^4(\Omega)}  \leq \frac{K\mu}{64}\|\be\|^2_\varepsilon+C\|\bu\|_2^2\|\bxi\|^2, \label{c21} \\
			|A_{7}|&\leq \|\bxi\|_{L^2(\Omega)}\|\|\nabla\be\|_{L^2(\Omega)}\|\bu\|_{L^\infty(\Omega)} \leq \frac{K\mu}{64}\|\be\|^2_\varepsilon+C\|\bu\|_2^2\|\bxi\|^2. \label{c22}
		\end{align}
		We next sum the last integral that is, $A_{10}$ over all $T$ and consider the contribution of this sum to one interior edge $e$. Assume that $e$ is shared by two triangles $T_i$ and $T_j$, with exterior normal $\bn_i$ and $\bn_j$. Then, we find
		\[\int_e (\bxi|_{T_i}\cdot \bn_i)\bu|_{T_i}\cdot \be|_{T_i}\,ds+\int_e(\bxi|_{T_j}\cdot \bn_j)\bu|_{T_j}\cdot \be|_{T_j}\,ds=\int_e[(\bxi\cdot \bn_e)\bu\cdot \be]\,ds. \]
		Thus, by using the trace inequality (\ref{T1}), we obtain
		\begin{align}
			A_{9}+A_{10}&=\frac{1}{2}\sum_{e\in\Gamma_h}\int_e\{\bxi\}\cdot \bn_e[\bu\cdot \be]\,ds \leq C \|\bu\|_{L^\infty(\Omega)}\sum_{e\in \Gamma_h}\frac{1}{|e|^{1/2}}\|[\be]\|_{L^2(e)}|e|^{1/2}\|\bxi\|_{L^2(e)}\nonumber\\
			& \leq C\|\bu\|_2 \sum_{T\in \mathcal{T}_h}\|\be\|_\varepsilon(\|\bxi\|_{L^2(T)}+h_T\|\nabla\bxi\|_{L^2(T)}) \leq  \frac{K\mu}{64}\|\be\|^2_\varepsilon+C \|\bu\|_2^2(\|\bxi\|^2+h^2\|\bxi\|_\varepsilon^2).\label{c23}
		\end{align}
		
		\noindent
		Incorporating (\ref{c4})-(\ref{c3}), (\ref{c11})-(\ref{c13}) and (\ref{c21})-(\ref{c23}) in (\ref{c0}), and thereby in (\ref{10.6}), and multiplying by $e^{2\al t}$ the resulting inequality, we observe that
		\begin{align}
			\frac{d}{dt}(e^{2\al t}\|\be\|^2)+(\mu K-2C(2)\al)e^{2\al t}\|\be\|_\varepsilon^2\leq & C(\|\bu\|_2^2+h^{-4}\|\bxi\|^2+h^{-2}\|\bxi\|_\varepsilon^2)e^{2\al t}\|\be\|^2\nonumber\\
			&+Ce^{2\al t}(\|\bu\|_2^2+h^{-2}\|\bxi\|_\varepsilon^2)(\|\bxi\|^2+h^2\|\bxi\|_\varepsilon^2).\label{c6}
		\end{align}
		Integrating (\ref{c6}) from $0$ to $t$, using
		$\|\be(0)\|\leq Ch^{2}$
		and  Gronwall's inequality, (\ref{13.2}) and Lemmas \ref{xierror} and \ref{exactpriori} in the resulting expression, we arrive at
		\begin{align*}
			e^{2\al t}\|\be\|^2+\mu K\int_0^te^{2\al s}\|\be\|_\varepsilon^2\,ds\leq Ce^{Ct}h^4.
		\end{align*}
		After multiplying the resulting inequality by $e^{-2\al t}$ and using the inverse relation (\ref{T6}), we obtain our desired estimate. This completes the proof.
	\end{proof}
	\noindent
 {\it Proof of Theorem \ref{semierror}.} Now, the proof of  Theorem \ref{semierror} follows by virtue of the Lemmas \ref{xierror} and \ref{etaerror}. \hfill{$\Box$}
	
	\noindent
	We would like to point out here that in the NIPG case, the estimates of $\be$ are suboptimal as they involve the estimates of $\bxi$ in both energy and $\bL^2$-norms, which are already shown as suboptimal in Remark \ref{Remark3}. Since $e=\bxi+\be$, we obtain the suboptimal estimates of semidiscrete velocity error in the NIPG case.
	\begin{Remark}
		The estimates of Theorem \ref{semierror} can be shown to be uniform (in time) under the smallness assumption on the data, that is,
		\begin{equation}
			N=\sup_{\bv_h,\bw_h\in\bV_h}\frac{c(\bw_h,\bv_h,\bw_h)}{\|\bw_h\|_\varepsilon^2\|\bv_h\|_\varepsilon}\quad \text{and}\quad \frac{N}{K^2\mu^2}\|\brf\|_{L^{2}(\Omega)}<1.\label{uniqueness}
		\end{equation}
	\end{Remark}
	\begin{proof}
		In order to derive estimates, which are valid uniformly for all $t>0$, let us rewrite the nonlinear terms in the following manner
		\begin{align*}
			c^{\bu_h}(\bu_h,\bu_h,\be)-c^{\bu_h}(\bu,\bu,\be)=- c^{\bu_h}(\bu,\bxi,\be)+c^{\bu_h}(\bxi,\bxi,\be)-c^{\bu_h}(\bxi,\bu,\be)-c^{\bu_h}(\be,\bu_h,\be)-c^{\bu_h}(\bv_h,\be,\be).
		\end{align*}
		From the proof of the Lemma \ref{etaerror}, we can derive the bounds as
		\[|c^{\bu_h}(\bu,\bxi,\be)+c^{\bu_h}(\bxi,\bu,\be)+c^{\bu_h}(\bxi,\bxi,\be)|\leq C(\|\bu\|_2+ch^{-1}\|\bxi\|_\varepsilon)(\|\bxi\|+h\|\bxi\|_\varepsilon)\|\be\|_\varepsilon,\]
		From the inequality (\ref{5.5}), we have
		\[c^{\bu_h}(\bv_h,\be,\be)\geq 0.\]
		Using the uniqueness condition, we find that
		\[c^{\bu_h}(\be,\bu_h,\be)\leq N\|\be\|_\varepsilon^2\|\bu_h\|_\varepsilon.\]
		We now modify the proof of Lemma \ref{etaerror} as follows:
		From (\ref{10.6}) and using Lemmas \ref{exactpriori} and \ref{xierror}, we obtain
		\begin{align}
			\frac{d}{dt}\|\be\|^2+2(\mu K-N\|\bu_h\|_\varepsilon)\|\be\|_\varepsilon^2\leq  Ch^2\|\be\|_\varepsilon. \label{uni1}
		\end{align}
		Multiply (\ref{uni1}) by $e^{2\al t}$ and integrate from $0$ to $t$. After a final multiplication of the resulting equation by $e^{-2\al t}$, we arrive at
		\begin{align*}
			\|\be\|^2+2e^{-2\al t}\int_0^te^{2\al s}(\mu K-N\|\bu_h\|_\varepsilon)\|\be(s)\|_\varepsilon^2\,ds\leq e^{-2\al t}\|\be(0)\|+2\al e^{-2\al t}\int_0^te^{2\al s}\|\be(s)\|^2\,ds\\
			+Ch^2e^{-2\al t}\int_0^te^{2\al s}\|\be\|_\varepsilon.
		\end{align*}
		Letting $t\rightarrow\infty$, applying the L'Hospital rule and using (\ref{uh4}), one can find
		\[\frac{1}{\al}\big(\mu K-\frac{N}{K\mu}\|\brf\|_{L^\infty(0,\infty;L^2(\Omega))}\big)\limsup_{t\rightarrow \infty}\|\be(t)\|_\varepsilon^2\leq \frac{Ch^2}{\al}\limsup_{t\rightarrow \infty}\|\be(t)\|_\varepsilon.\]
		Due to the uniqueness condition (\ref{uniqueness}), there holds
		\[\limsup_{t\rightarrow \infty}\|\be(t)\|_\varepsilon\leq Ch^2.\]
		Therefore,
		\[\limsup_{t\rightarrow \infty}\|\be(t)\|\leq Ch^2.\]
		Together with the estimate of $\bxi$ from Lemma \ref{xierror}, we find
		\[\limsup_{t\rightarrow \infty}\|\bu(t)-\bu_h(t)\|\leq Ch^2.\]
		Here, the constant $C$ is valid uniformly for all $t>0$.
	\end{proof}
	
	\section{Error estimates for pressure}\label{s5}
	\se
	In this section, we derive error estimates for the semi-discrete discontinuous Galerkin approximation of the pressure. Before proving our main theorem we need some auxiliary lemmas.   
	\begin{lemma}\label{etL2L2err}
		The velocity error $\de=\bu-\bu_h$ satisfies, for $0<t<T$,
		\begin{equation}
			\int_0^te^{2\al s}\|\de_t(s)\|^2\,ds\leq Ce^{Ct}h^2.\label{etL2L2}
		\end{equation}
	\end{lemma}
	\begin{proof}
		Let us denote $\bchi=\bS_h\bu-\bu_h$. From the equations for $\bu$, $\bu_h$ and $\bS_h\bu$, that is, (\ref{8.5}), (\ref{u_h}) and (\ref{9.3}), respectively, and for $\bphi_h\in\bV_h$, we obtain
		\begin{align*}
			(\bchi_t,\bphi_h)+\mu(a(\bchi,\bphi_h)+J_0(\bchi,\bphi_h))+c^{\bu}(\bu,\bu,\bphi_h)-c^{\bu_h}(\bu_h,\bu_h,\bphi_h)=-(\bzeta_t,\bphi_h).
		\end{align*}
		Choose $\bphi_h=\bchi_t$ in the above equality to obtain
		\begin{align}
			\|\bchi_t\|^2+\frac{\mu}{2}\frac{d}{dt}\|\bchi\|_\varepsilon^2+c^{\bu}(\bu,\bu,\bchi_t)-c^{\bu_h}(\bu_h,\bu_h,\bchi_t)=-(\bzeta_t,\bchi_t).\label{et1}
		\end{align}
		We can drop the superscripts from $c(\cdot,\cdot,\cdot)$ and rewrite the nonlinear terms as
		\[c^{\bu}(\bu,\bu,\bchi_t)-c^{\bu_h}(\bu_h,\bu_h,\bchi_t)=-c(\de,\de,\bchi_t)+c(\de,\bu,\bchi_t)+c(\bu,\de,\bchi_t).\]
		By using $L^p$ bound, Lemmas \ref{trace} and \ref{distrace}, and Theorem \ref{semierror}, we  bound $c(\de,\de,\bchi_t)$ similar to Lemma \ref{etaerror} as follows
		\begin{align}
			|c(\de,\de,\bchi_t)|\leq C\|\de\|_\varepsilon\|\bchi_t\|\leq \frac{1}{6}\|\bchi_t\|^2+C\|\de\|_\varepsilon^2.\label{pr1}
		\end{align}
		Since $\bu$ is continuous, Lemma \ref{trace}, Sobolev's inequality and (\ref{priori1}) yield
		\begin{align}
			|c(\de,\bu,\bchi_t)|=&\bigg|\sum_{T\in\mathcal{T}_h}\int_T(\de\cdot \nabla \bu)\cdot \bchi_t\,dT+\frac{1}{2}\sum_{T\in\mathcal{T}_h}\int_T(\nabla\cdot \de)\bu\cdot \bchi_t\,dT-\frac{1}{2}\sum_{e\in\Gamma_h}\int_e[\de]\cdot n_e\{\bu\cdot\bchi_t\}\,ds\bigg|\nonumber\\
			\leq &\sum_{T\in\mathcal{T}_h}\|\de\|_{L^4(T)}\|\nabla\bu\|_{L^4(T)}\|\bchi_t\|_{L^2(T)}+C\|\bu\|_{L^\infty(\Omega)}\sum_{T\in\mathcal{T}_h}\|\nabla\de\|_{L^2(T)}\|\bchi_t\|_{L^2(T)}\nonumber\\
			&+C\|\bu\|_{L^\infty(\Omega)}\sum_{T\in\mathcal{T}_h}\frac{1}{|e|^{1/2}}\|[\de]\|_{L^2(e)}\|\bchi_t\|_{L^2(T)}\leq \frac{1}{6}\|\bchi_t\|^2+C\|\de\|_\varepsilon^2.\label{pr2}
		\end{align}
		Similarly, we can  bound
		\begin{align}
			|c(\bu,\de,\bchi_t)|\leq \frac{1}{6}\|\bchi_t\|^2+C\|\de\|_\varepsilon^2.\label{pr3}
		\end{align}
		Apply (\ref{pr1})--(\ref{pr3}) in (\ref{et1}) and multiply the resulting inequality by $e^{2\al t}$. Then, integrating from $0$ to $t$ with respect to time, we obtain
		\begin{align}
			\int_0^t e^{2\al s}\|\bchi_t\|^2\,ds+\mu e^{2\al t}\|\bchi\|_\varepsilon^2\leq \mu\|\bchi(0)\|_\varepsilon^2+2\al \int_0^t e^{2\al s}\|\bchi\|_\varepsilon^2\,ds+C\int_0^t e^{2\al s}\|\bzeta_t\|^2\,ds+C\int_0^t e^{2\al s}\|\de\|_\varepsilon^2\,ds.\label{et3}
		\end{align}
		Again, by using the estimates (\ref{9.4}) and (\ref{priori2}), we have
		\[\int_0^t e^{2\al s}\|\bchi\|_\varepsilon^2\,ds\leq Ce^{2\al t}h^2+\int_0^t e^{2\al s}\|\de\|_\varepsilon^2\,ds.\]
		Using (\ref{9.5}), (\ref{priori3}) and Theorem \ref{semierror} in (\ref{et3}), we obtain
		\[\int_0^t e^{2\al s}\|\bchi_t\|^2\,ds\leq Ce^{Ct}h^2.\]
		Furthermore, a use of triangle inequality, estimate (\ref{9.5}) and (\ref{priori3}) lead to
		\[\int_0^t e^{2\al s}\|\de_t\|^2\,ds\leq Ce^{Ct}h^2.\]
	\end{proof}
	\begin{lemma}\label{etLinfL2}
		The error $\de=\bu-\bu_h$ in approximating the velocity satisfies for $t>0$
		\begin{equation*}
			\|\de_t(t)\|\leq Ce^{Ct}h(\tau(t))^{-1/2}, 
		\end{equation*}
		where $\tau(t)=\min\{t,1\}$.
	\end{lemma}
	\begin{proof}
		The error equation in $\de$ obtained from (\ref{8.5}) and (\ref{u_h}) is
		\begin{align}
			(\de_t,\bphi_h)+\mu(a(\de,\bphi_h)+J_0(\de,\bphi_h))+c^{\bu}(\bu,\bu,\bphi_h)-c^{\bu_h}(\bu_h,\bu_h,\bphi_h)+b(\bphi_h,p)=0\quad \forall \bphi_h\in \bV_h.\label{et4}
		\end{align}
		Since $\bu$ is continuous, we can drop the superscripts of the nonlinear terms. 
		Differentiate (\ref{et4}) with respect to $t$ and choose $\bphi_h=\bP_h\de_t$. Then using the definition of $\bP_h$ and Lemma \ref{coer}, we find
		\begin{align}
			\frac{d}{dt}\|\de_t\|^2+2\mu K\|\bP_h\de_t\|_\varepsilon^2 &\leq \frac{d}{dt}\|\bu_t-\bP_h\bu_t\|^2+2\mu(a(\bu_t-\bP_h\bu_t,\bP_h\de_t)+J_0(\bu_t-\bP_h\bu_t,\bP_h\de_t))\nonumber\\
			&+2(c(\bu_{ht},\bu_h,\bP_h\de_t)+c(\bu_h,\bu_{ht},\bP_h\de_t)-c(\bu_t,\bu,\bP_h\de_t) -c(\bu,\bu_t,\bP_h\de_t))-2b(\bP_h\de_t,p_t).\label{et5}
		\end{align}
		We rewrite the nonlinear terms as
		\begin{align}
			c(\bu_{ht},\bu_h,\bP_h\de_t)+(c(\bu_h,\bu_{ht},\bP_h\de_t)-c(\bu_t,\bu,\bP_h\de_t)-c(\bu,\bu_t,\bP_h\de_t)=-c(\bu_{ht},\de,\bP_h\de_t)-c(\de_t,\bu,\bP_h\de_t)\nonumber\\
			-c(\de,\bu_{ht},\bP_h\de_t)-c(\bu,\de_t,\bP_h\de_t).\label{pr4}
		\end{align}
		Using Cauchy-Schwarz's inequality,  Young's inequality and Lemmas \ref{trace} and \ref{distrace}, the nonlinear terms on the right hand side of (\ref{pr4}) can be bounded as in Lemma \ref{etaerror}. Therefore, we have
		\begin{align}
			|c(\bu_{ht},\de,\bP_h\de_t)| &\leq  \frac{\mu K}{12}\|\bP_h\de_t\|_\varepsilon^2+Ch^{-1}(\|\de\|^2+h^2\|\de\|_\varepsilon^2)\|\bu_{ht}\|_\varepsilon^2,\\
			|c(\de_t,\bu,\bP_h\de_t)|&\leq \frac{\mu K}{12}\|\bP_h\de_t\|_\varepsilon^2+C(\|\de_t\|^2+h^2\|\bu_t\|_2^2)\|\bu\|_2^2,\\
			|c(\bu,\de_t,\bP_h\de_t)|&\leq  \frac{\mu K}{12}\|\bP_h\de_t\|_\varepsilon^2+C(\|\de_t\|^2+h^2\|\bu_t\|_2^2)\|\bu\|_2^2,\\
			|c(\de,\bu_{ht},\bP_h\de_t)| &\leq  \frac{\mu K}{12}\|\bP_h\de_t\|_\varepsilon^2+Ch^{-1}(\|\de\|^2+h^2\|\de\|_\varepsilon^2)\|\bu_{ht}\|_\varepsilon^2.
		\end{align}
		The other terms on the right hand side of (\ref{et5}) can be bounded as in Lemma \ref{xiL2err}:
		\begin{align}
			2\mu|a(\bu_t-\bP_h\bu_t,\bP_h\de_t)+J_0(\bu_t-\bP_h\bu_t,\bP_h\de_t)|\leq \frac{K\mu}{12}\|\bP_h\de_t\|_\varepsilon^2+Ch^{2}|\bu_t|_2^2,\\ 
			2|b(\bP_h\de_t,p_t))|=2|b(\bP_h\de_t,p_t-r_h(p_t)))|\leq\frac{K\mu}{12}\|\bP_h\de_t\|_\varepsilon^2+Ch^{2}|p_t|_1^2.\label{pr6}
		\end{align}
		Using the bounds from (\ref{pr4})-(\ref{pr6}) in (\ref{et5}), we obtain
		\begin{align}
			\frac{d}{dt}\|\de_t\|^2+\mu K\|\bP_h\de_t\|_\varepsilon^2\leq &\frac{d}{dt}\|\bu_t-\bP_h\bu_t\|^2+C\|\de_t\|^2+Ch^{-1}(\|\de\|^2+h^2\|\de\|_\varepsilon^2)\|\bu_{ht}\|_\varepsilon^2+Ch^{2}(\|\bu_t\|_2^2+\|p_t\|_1^2).\label{et6}
		\end{align}
		Multiply (\ref{et6}) by $\sigma(t)$ and integrate with respect to time, to write
		\begin{align*}
			\sigma(t)\|\de_t\|^2+\mu K\int_0^t\sigma(s)\|\bP_h\de_t\|_\varepsilon^2\,ds\leq \int_0^t(\sigma_s(s)+\sigma(s))\|\de_t\|^2\,ds+\sigma(t)\|\bu_t-\bP_h\bu_t\|^2\\
			+Ch^{-1}\int_0^t\sigma(s)(\|\de\|^2+h^2\|\de\|_\varepsilon^2)\|\bu_{ht}\|_\varepsilon^2\,ds+Ch^{2}\int_0^t\sigma(s)(|\bu_t|_2^2+|p_t|_1^2)\,ds.
		\end{align*}
		Finally, from Lemmas \ref{etL2L2err} and \ref{P_h}, Theorem \ref{semierror}, (\ref{uh3}) and (\ref{priori3}), we conclude the rest of the proof.
	\end{proof}
	\begin{theorem}\label{semipressure}
		Under the assumptions of Theorem \ref{semierror} , there exists a positive constant $C$, independent of $h$, such that the following error estimate holds:
		\begin{align*}
			\|(p-p_h)(t)\|\leq Ce^{Ct}ht^{-1/2}, \quad 0<t<T.
		\end{align*}
	\end{theorem}
	\begin{proof}
		From (\ref{8.5}), (\ref{8.6}), (\ref{8.8}) and (\ref{8.9}), we can write the error equation as follows:
		\begin{align}
			(\bu_{ht}-\bu_t,\bv_h)+\mu(a(\bu_h-\bu,\bv_h)+J_0(\bu_h-\bu,\bv_h))+c^{\bu_h}(\bu_h,\bu_h,\bv_h)-c^{\bu}(\bu,\bu,\bv_h)\nonumber\\
			-b(\bv_h,p-r_h(p))=-b(\bv_h,p_h-r_h(p)), \quad \forall \bv_h\in\bX_h.\label{30.16}
		\end{align}
		By virtue of the inf-sup condition in Lemma \ref{inf-sup}, there exists $\bv_h\in \bX_h$ such that
		\begin{equation}
			b(\bv_h,p_h-r_h(p))=-\|p_h-r_h(p)\|^2,\quad \|v_h\|_\varepsilon\leq \frac{1}{\beta_0}\|p_h-r_h(p)\|.\label{30.17}
		\end{equation}
		Therefore, from (\ref{30.16}), we obtain
		\begin{align}
			\|p_h-r_h(p)\|^2=&(\bu_{ht}-\bu_t,\bv_h)+\mu(a(\bu_h-\bu,\bv_h)+J_0(\bu_h-\bu,\bv_h))\nonumber\\
			&+c(\bu_h,\bu_h,\bv_h)-c(\bu,\bu,\bv_h)-b(\bv_h,p-r_h(p)).\label{30.18}
		\end{align}
		The terms on the right hand side of (\ref{30.18}) can be bounded as in Lemmas \ref{xiL2err} and \ref{etaerror}. Then, the inequality (\ref{30.18}) becomes
		\begin{align*}
			\|p_h-r_h(p)\|^2\leq C\|\bu_{ht}-\bu_t\|^2+C\|\bu_h-\bu\|_\varepsilon^2+Ch^2(|\bu|_2^2+|p|_1^2)+C\|\bu_h-\bu\|^2.
		\end{align*}
		Using the triangle inequality, we obtain
		\begin{align*}
			\|p-p_h\|^2\leq C\|\bu_{ht}-\bu_t\|^2+C\|\bu_h-\bu\|_\varepsilon^2+Ch^2(|\bu|_2^2+|p|_1^2)+C\|\bu_h-\bu\|^2.
		\end{align*}
		Combining Lemma \ref{etLinfL2}, Theorem \ref{semierror} and  (\ref{priori1}) with the above inequality we obtain our desired pressure error estimate.
	\end{proof}

 It can be noted that the pressure estimates in the NIPG case will be suboptimal due to their dependence on velocity error estimates, which are already discussed as suboptimal in Section \ref{s4}.
	\section{Fully discrete approximation and error estimates}\label{s6}
	\se
	In this section, we analyse the backward Euler method for temporal discretization of the semidiscrete approximations to the discontinuous Galerkin time dependent Navier-Stokes equations. For time discretization,  let $\Delta t$, $0<\Delta t<1$, denote the time step, $t_n=n\Delta t$, $n\geq 0$, and $0=t_0<t_1<\cdots<t_M=T$ be a subdivision of the interval $(0,T)$. We denote the function $\phi(t)$ evaluated at time $t_n$ by $\phi^n$. We define for a sequence $\{\phi^n\}_{n\geq 0}\subset \bV_h$,
	\[\partial_t\phi^n=\frac{1}{\Delta t}(\phi^n-\phi^{n-1}).\]
	We describe below the backward Euler scheme for the semi-discrete problem (\ref{8.8})-(\ref{8.9}) as follows:\\
	Given $\bU^0$, find $\{\bU^n\}_{n\geq 1}\in \bX_h$ and $\{P^n\}_{n\geq 1}\in M_h$ such that
	\begin{align}
		(\partial_t \bU^n,\bv_h)+\mu(a(\bU^n,\bv_h)+J_0(\bU^n,\bv_h))+c^{\bU^{n-1}}(\bU^{n-1},\bU^n,\bv_h) +b(\bv_h,P^n)&=(\brf^n,\bv_h),\quad \forall \bv_h\in \bX_h,\label{fullyxh1}\\
		b(\bU^n,q_h)&=0,\quad \forall q_h\in M_h. \label{fullyxh2}
	\end{align}
	Note that $\bU^0=\bu_h(0)=\bP_h\bu_0$.\\
	Now, for $\bv_h\in \bV_h$, we seek $\{\bU^n\}_{n\geq 1}\in \bV_h$ such that
	\begin{align}
		(\partial_t \bU^n,\bv_h)+\mu(a(\bU^n,\bv_h)+J_0(\bU^n,\bv_h &))+c^{\bU^{n-1}}(\bU^{n-1},\bU^n,\bv_h) =(\brf^n,\bv_h),\quad \forall \bv_h\in \bV_h. \label{fullysol}
	\end{align}
	\begin{lemma}\label{fullypriorisol}
		The solution $\{\bU^n\}_{n\geq 1}$ of (\ref{fullysol}) satisfies the following estimates:
		\begin{align*}
			\|\bU^n\|^2+e^{-2\al t_M}\Delta t\sum_{n=1}^{M}e^{2\al t_n}\|\bU^n\|^2_\varepsilon &\leq C,\quad n=0,1,...,M,
		\end{align*}
		where $C$ depends on the given data.
	\end{lemma}
	The above estimates of $\bU^n$ can be easily derived by choosing $\bv_h=\bU^n$ in (\ref{fullysol}) and using Lemma \ref{coer}.
	Next, we discuss the error estimate of the backward Euler method. For the error analysis, we set, for fixed $n\in\mathbb{N},\, 1<n\leq M,\, \de_n=\bU^n-\bu_h(t_n)=\bU^n-\bu_h^n$.
	Considering the semi-discrete formulation (\ref{u_h}) at $t=t_n$ and subtracting from (\ref{fullysol}), we arrive at
	\begin{align}
		(\partial_t\de_n,\bv_h)+\mu( & a(\de_n,\bv_h)+J_0(\de_n,\bv_h))=(\bu_{ht}^n,\bv_h) -(\partial_t\bu_h^n,\bv_h)+c^{\bu_h^n}(\bu_h^n,\bu_h^n,\bv_h)-c^{\bU^{n-1}}(\bU^{n-1},\bU^n,\bv_h).\label{11.0}
	\end{align}
	Using Taylor's series expansion, we observe that 
	\begin{equation}\label{11.1}
		(\bu^n_{ht},\bv_h)-(\partial_t\bu_h^n,\bv_h)=\frac{1}{\Delta t}\int_{t_{n-1}}^{t_n} (t-t_{n-1})(\bu_{htt},\bv_h)\, dt.
	\end{equation}
	\begin{lemma}\label{fully}
		Under the assumptions of Theorem \ref{semierror}, there is a positive constant $C$, independent of $h$ and $\Delta t$, such that
		\begin{align*}
			\|\de_n\|+\left(\mu K e^{-2\al t_M}\Delta t\sum_{n=1}^{M}e^{2\al t_n}\|\de_n\|^2_\varepsilon\right)^{1/2}\leq Ce^{CT}\Delta t.
		\end{align*}
	\end{lemma}
	\begin{proof}
		We put $\bv_h=\de_n$ in error equation (\ref{11.0}). With the observation 
		\[(\partial _t\de_n,\de_n)=\frac{1}{\Delta t}(\de_n-\de_{n-1},\de_n)\geq \frac{1}{2\Delta t}(\|\de_n\|^2-\|\de_{n-1}\|^2)=\frac{1}{2}\partial_t \|\de_n\|^2,\]
		and using Lemma \ref{coer}, we find that
		\begin{align}
			\partial_t\|\de_n\|^2+2K\mu\|\de_n\|^2_\varepsilon\leq &2(\bu^n_{ht},\de_n)-2(\partial_t\bu_h^n,\de_n)  + 2c^{\bu_h^n}(\bu_h^n,\bu_h^n,\de_n)-2c^{\bU^{n-1}}(\bU^{n-1},\bU^n,\de_n).\label{11.2}
		\end{align}
		The nonlinear terms in (\ref{11.2}) can be written as
		\begin{align}
			c^{\bu_h^n}(\bu_h^n,\bu_h^n,\de_n)-c^{\bU^{n-1}}(\bU^{n-1},\bU^n,\de_n)=&c^{\bu_h^n}(\bu_h^n,\bu_h^n,\de_n)-c^{\bu_h^n}(\bU^{n-1},\bU^n,\de_n)\nonumber\\
			&+c^{\bu_h^n}(\bU^{n-1},\bU^n,\de_n)-c^{\bU^{n-1}}(\bU^{n-1},\bU^n,\de_n).\label{en1}
		\end{align}
		We now drop the superscripts for the first two nonlinear terms on the right hand side of (\ref{en1}) and rewrite them as
		\begin{align}
			c(\bu_h^n, \bu_h^n,\de_n)-c(\bU^{n-1}  ,\bU^n,\de_n)= &-c(\bU^{n-1},\de_n,\de_n) +c(\de_{n-1},\bu^n-\bu_h^n,\de_n)-c(\de_{n-1},\bu^n,\de_n)\nonumber\\
			&+c(\bu_h^n-\bu_h^{n-1},\bu_h^n-\bu^n,\de_n)+c(\bu_h^n-\bu_h^{n-1},\bu^n,\de_n).\label{en2}
		\end{align}
		From (\ref{5.5}), we have $c(\bU^{n-1},\de_n,\de_n)\geq 0$. The term $c(\de_{n-1},\bu^n-\bu_h^n,\de_n)$ in (\ref{en2}) can be bounded exactly like the term $c(\be,\bxi,\be)$ in Lemma \ref{etaerror}. Cauchy-Schwarz's inequality and Theorem \ref{semierror} yield
		\begin{align}
			|c(\de_{n-1},\bu^n-\bu_h^n,\de_n)|\leq C\|\de_{n-1}\|\|\de_n\|_\varepsilon\leq \frac{K\mu}{64}\|\de_n\|^2_\varepsilon+\frac{C}{\mu}\|\de_{n-1}\|^2.\label{en3}
		\end{align}
		Using estimate (\ref{tri1}) and Cauchy-Schwarz's inequality and Sobolev's inequality yield
		\begin{align}
			|c(\de_{n-1},\bu^n,\de_n)|\leq C\|\de_{n-1}\|_{L^2(\Omega)}|\bu^n|_{W^{1,4}(\Omega)}\|\de_n\|_\varepsilon\leq \frac{K\mu}{64}\|\de_n\|^2_\varepsilon+C\|\de_{n-1}\|^2. \label{en4}
		\end{align}
		Next, we are going to use the fact that, $\bu_h^n-\bu_h^{n-1}=\displaystyle\int_{t_{n-1}}^{t_n}\bu_{ht}(t)\,dt$. Using H{\"o}lder's inequality, Young's inequality,  Lemmas \ref{trace} and \ref{distrace}, (\ref{Lp}) and Theorem \ref{semierror},  the fourth term on the right hand side in (\ref{en2}) is bounded, following the form of $c(\cdot,\cdot,\cdot)$ in (\ref{4.8}) as  
		\begin{align}
			|c(\bu_h^n-\bu_h^{n-1},\bu_h^n-\bu^n,\de_n)|&\leq  \int_{t_{n-1}}^{t_n}\sum_{T\in \mathcal{T}_h}\|\bu_{ht}\|_{L^4(T)}\|\nabla(\bu_h^n-\bu^n)\|_{L^2(T)}\|\de_n\|_{L^4(T)}\,dt\nonumber\\
			+\int_{t_{n-1}}^{t_n}\sum_{e\in\Gamma_h}\|\bu_{ht}\|_{L^4(e)}\|[\bu_h^n &-\bu^n]\|_{L^2(e)}\|\de_n\|_{L^4(e)}\,dt +C\int_{t_{n-1}}^{t_n}\sum_{T\in \mathcal{T}_h}\|\nabla\bu_{ht}\|_{L^2(T)}\|\bu_h^n-\bu^n\|_{L^4(T)}\|\de_n\|_{L^4(T)}\,dt\nonumber\\
			+C\int_{t_{n-1}}^{t_n} \sum_{e\in \Gamma_h}\|[\bu_{ht}]\|_{L^2(e)}&\|\bu_h^n-\bu^n\|_{L^4(e)}\|\de_n\|_{L^4(e)}\,dt \nonumber\\
			\leq  \frac{K\mu}{64}\|\de_n\|^2_\varepsilon +\frac{C}{\mu}\Delta t &\|\bu_{ht}\|^2_{L^2(t_{n-1},t_n;\varepsilon)}.
		\end{align}
		As $\bu^n$ has zero jump, the last nonlinear term $c(\bu_h^n -\bu_h^{n-1},\bu^n,\de_n)$ on the right hand side in (\ref{en2}) can be bounded, following the form of $c(\cdot,\cdot,\cdot)$ in (\ref{4.8}) and  using (\ref{Lp}) and Lemma \ref{distrace} as
		\begin{align}
			|c (\bu_h^n &-\bu_h^{n-1} ,\bu^n,\de_n)|\nonumber\\
			\leq &\Delta t^{1/2}\|\bu_{ht}\|_{L^2(t_{n-1},t_n;L^2(\Omega))}\|\nabla\bu^n\|_{L^4(\Omega)}\|\de_n\|_{L^4(\Omega)}+C\Delta t^{1/2} \|\bu_{ht}\|_{L^2(t_{n-1},t_n;\varepsilon)}\|\bu^n\|_{L^\infty([0,T]\times\Omega)}\|\de_n\|_{L^2(\Omega)}\nonumber\\
			& +C\int_{t_{n-1}}^{t_n}\sum_{T\in \mathcal{T}_h}\|\bu^n\|_{L^\infty([0,T]\times\Omega)}|e|^{-1/2}\|[\bu_{ht}]\|_{L^2(e)}|e|^{1/2}h^{-1/2}\|\de_n\|_{L^2(T)}\,dt\nonumber\\
			\leq &\frac{K\mu}{64}\|\de_n\|^2_\varepsilon+\frac{C}{\mu}\Delta t\|\bu_{ht}\|^2_{L^2(t_{n-1},t_n;\varepsilon)}.
		\end{align}
		From Proposition 4.10 in \cite{GR09}, we can easily see that
		\begin{align*}
			|c^{\bu_h^n}(\bU^{n-1},\bU^n,\de_n)-c^{\bU^{n-1}}(\bU^{n-1},\bU^n,\de_n)|\leq \|\bu_h^n-\bU^{n-1}\|_{L^4(\Omega)}\|\bU^n\|_\varepsilon\|\de_n\|_\varepsilon.
		\end{align*}
		A use of the triangle inequality, Theorem 3.8 in \cite{GRW05}, estimate (\ref{uh1}), Lemma \ref{fullypriorisol},  and Cauchy-Schwarz's inequality yield
		\begingroup
		\allowdisplaybreaks
		\begin{align}
			| c^{\bu_h^n}&(\bU^{n-1},\bU^n,\de_n)-c^{\bU^{n-1}}(\bU^{n-1},\bU^n,\de_n)|\nonumber\\
			\leq &\|\bu_h^n-\bu_h^{n-1}\|_{L^4(\Omega)}\|\de_n\|_\varepsilon\|\de_n\|_{L^4(\Omega)} +\|\de_{n-1}\|_{L^4(\Omega)}\|\de_n\|_\varepsilon\|\de_n\|_{L^4(\Omega)}\nonumber\\
			&+ \|\bu_h^n-\bu_h^{n-1}\|_{L^4(\Omega)}\|\bu_h^n\|_\varepsilon\|\de_n\|_{L^4(\Omega)} +\|\de_{n-1}\|_{L^4(\Omega)}\|\bu_h^n\|_\varepsilon\|\de_n\|_{L^4(\Omega)}\nonumber\\
			\leq &\|\bu_h^n-\bu_h^{n-1}\|^{1/2}\|\bu_h^n-\bu_h^{n-1}\|_\varepsilon^{1/2}\|\de_n\|_\varepsilon^{3/2}\|\de_n\|^{1/2} +\|\de_{n-1}\|^{1/2}\|\de_{n-1}\|_\varepsilon^{1/2}\|\de_n\|_\varepsilon^{3/2}\|\de_n\|^{1/2}\nonumber\\
			&+ C\|\bu_h^n-\bu_h^{n-1}\|_\varepsilon\|\de_n\|_\varepsilon +C\|\de_{n-1}\|^{1/2}\|\de_{n-1}\|_\varepsilon^{1/2}\|\de_n\|_\varepsilon\nonumber\\
			\leq &  \frac{K\mu}{64}\|\de_n\|^2_\varepsilon+\frac{K\mu}{64}\|\de_{n-1}\|^2_\varepsilon+C\|\de_{n-1}\|^2+C\|\de_{n-1}\|^2\|\de_{n-1}\|^2_\varepsilon+C\Delta t\|\bu_{ht}\|^2_{L^2(t_{n-1},t_n;\varepsilon)}.
		\end{align}
		From (\ref{11.1}), we have 
		\begin{align}
			2(\bu^n_{ht},\de_n)-2(\partial_t\bu_h^n,\de_n)& \leq C\Delta t^{1/2}\left(\int_{t_{n-1}}^{t_n}\|\bu_{htt}(t)\|^2_{-1,h}\,dt\right)^{1/2}\|\de_n\|_\varepsilon\nonumber\\
			&\leq  \frac{K\mu}{64}\|\de_n\|^2_\varepsilon+\frac{C}{\mu}\Delta t\int_{t_{n-1}}^{t_n}\|\bu_{htt}(t)\|^2_{-1,h}\,dt.\label{en9}
		\end{align}
		\endgroup
		Combine (\ref{en1})-(\ref{en9}), multiply (\ref{11.2}) by $\Delta t e^{2 \al n \Delta t}$ and sum over $1\leq n\leq m\leq M$, where $T=M\Delta t$ and observe that
		\begin{align*}
			\sum_{n=1}^m \Delta t e^{2\al n \Delta t}\partial_t\|\de_n\|^2=e^{2\al m \Delta t}\|\de_m\|^2-\sum_{n=1}^{m-1}  e^{2\al n \Delta t}( e^{2\al \Delta t}-1)\|\de_m\|^2
		\end{align*}
		to obtain
		\begin{align}
			e^{2\al m \Delta t}\|\de_m\|^2+K\mu\Delta t\sum_{n=1}^me^{2\al n\Delta t}\|\de_n\|_\varepsilon^2\leq \sum_{n=1}^{m-1}  e^{2\al n \Delta t}( e^{2\al \Delta t}-1)\|\de_n\|^2+\Delta t\sum_{n=1}^m e^{2\al n \Delta t}(C+\|\de_{n-1}\|_\varepsilon^2)\|\de_{n-1}\|^2\nonumber\\
			+ C{\Delta t}^2\sum_{n=1}^m e^{2\al n \Delta t}\int_{t_{n-1}}^{t_n}\|\bu_{ht}(t)\|_\varepsilon^2\,dt+C{\Delta t}^2\sum_{n=1}^m e^{2\al n \Delta t}\int_{t_{n-1}}^{t_n}\|\bu_{htt}(t)\|_{-1,h}^2\,dt.\label{11.5}
		\end{align}
		Note that as $e^{2\al \Delta t}-1\leq C(\al)\Delta t$, so the first term on right hand side of (\ref{11.5}) can be merged with the second term on right hand side of (\ref{11.5}). We bound the terms involving $\bu_h$ using Lemma \ref{uhpriori}. Observe that 
		\begin{align*}
			\sum_{n=1}^m e^{2\al n \Delta t}\int_{t_{n-1}}^{t_n}\|\bu_{ht}(t)\|_\varepsilon^2\,dt=\sum_{n=1}^m\int_{t_{n-1}}^{t_n}e^{2\al (t_n-t)}e^{2\al t}\|\bu_{ht}(t)\|_\varepsilon^2\,dt\leq e^{2\al \Delta t}\int_0^{t_m}e^{2\al t}\|\bu_{ht}(t)\|_\varepsilon^2\,dt\leq Ce^{2\al (m+1)\Delta t}.
		\end{align*}
		Similarly, the last term on the right hand side of (\ref{11.5}) can be bounded by $C{\Delta t}^2e^{2\al (m+1)\Delta t}$. Therefore, from (\ref{11.5}), we obtain 
		\[e^{2\al m \Delta t}\|\de_m\|^2+K\mu\Delta t\sum_{n=1}^me^{2\al n\Delta t}\|\de_n\|_\varepsilon^2\leq C\Delta t\sum_{n=1}^m e^{2\al (n-1) \Delta t}(C+\|\de_{n-1}\|_\varepsilon^2)\|\de_{n-1}\|^2+C{\Delta t}^2e^{2\al (m+1)\Delta t}.\]
		A use of the  discrete Gronwall's lemma completes the rest of the proof.
	\end{proof}
	\begin{theorem}\label{velofinal}
		Under the assumptions of the Theorem \ref{semierror} and Lemma \ref{fully}, the following holds true:
		\[\|\bu(t_n)-\bU^n\|\leq Ce^{CT}(h^2+\Delta t),\]
		where $C$ is a positive constant, independent of $h$ and $\Delta t$.
	\end{theorem}
	\begin{proof}
		Combine  Theorem \ref{semierror} with Lemma \ref{fully} to complete the proof.
	\end{proof}
	\begin{lemma}\label{L2L2delt}
		Under the assumptions of Lemma \ref{fully}, the error $\de_n=\bU^n-\bu_h^n$, satisfies
		\begin{align*}
			e^{-2\al t_m}\Delta t\sum_{n=1}^m  e^{2\al t_n}\|\partial_t\de_n\|^2+\mu K\|\de_m\|_\varepsilon^2\leq C\Delta t.
		\end{align*}
	\end{lemma}
	\begin{proof}
		In (\ref{11.0}), choose $\bv_h=\partial_t\de_n$ to obtain
		\begin{align}
			\|\partial_t\de_n\|^2+\mu( & a(\de_n,\partial_t\de_n)+J_0(\de_n,\partial_t\de_n))=(\bu_{ht}^n,\partial_t\de_n)\nonumber\\
			&-(\partial_t\bu_h^n,\partial_t\de_n)+c^{\bu_h^n}(\bu_h^n,\bu_h^n,\partial_t\de_n)-c^{\bU^{n-1}}(\bU^{n-1},\bU^n,\partial_t\de_n).\label{ent1}
		\end{align}
		The nonlinear terms in (\ref{ent1}) can be written as
		\begin{align}
			c^{\bu_h^n}(\bu_h^n,\bu_h^n,\partial_t\de_n)-c^{\bU^{n-1}}(\bU^{n-1},\bU^n,\partial_t\de_n)=&c^{\bu_h^n}(\bu_h^n,\bu_h^n,\partial_t\de_n)-c^{\bu_h^n}(\bU^{n-1},\bU^n,\partial_t\de_n)\nonumber\\
			&+c^{\bu_h^n}(\bU^{n-1},\bU^n,\partial_t\de_n)-c^{\bU^{n-1}}(\bU^{n-1},\bU^n,\partial_t\de_n).\label{ent2}
		\end{align}
		We now drop the superscripts for the first two nonlinear terms on the right hand side of (\ref{ent2}) and rewrite them as
		\begin{align}
			c(\bu_h^n, \bu_h^n,\partial_t\de_n)-c(\bU^{n-1} ,\bU^n,\partial_t\de_n)= &-c(\bu_h^{n-1}-\bu^{n-1},\de_n,\partial_t\de_n) -c(\bu^{n-1},\de_n,\partial_t\de_n)\nonumber\\
			&+c(\bu_h^n-\bu_h^{n-1},\bu_h^n-\bu^n,\partial_t\de_n)+c(\bu_h^n-\bu_h^{n-1},\bu^n,\partial_t\de_n)\nonumber\\
			&-c(\de_{n-1},\bu_h^n-\bu^n,\partial_t\de_n)-c(\de_{n-1},\bu^n,\partial_t\de_n)-c(\de_{n-1},\de_n,\partial_t\de_n).\label{ent3}
		\end{align}
		$L^p$ bound (\ref{Lp}), Lemma \ref{trace}, Lemma \ref{distrace}, Theorem \ref{semierror} and Sobolev's inequalities give the bounds for the nonlinear terms in the right hand side of (\ref{ent3}) except the last term,  as in Lemma \ref{fully}. 
		\begin{align}
			|c (\bu_h^{n-1}-\bu^{n-1},\de_n,\partial_t\de_n)|&\leq \frac{1}{64}\|\partial_t\de_n\|^2+C \|\de_n\|_\varepsilon^2,\label{ent5}\\
			|c(\bu^{n-1},\de_n,\partial_t\de_n)|&\leq \frac{1}{64}\|\partial_t\de_n\|^2+C \|\de_n\|_\varepsilon^2,\\
			|c(\bu_h^n-\bu_h^{n-1},\bu_h^n-\bu^n,\partial_t\de_n)+c(\bu_h^n-\bu_h^{n-1},\bu^n,\partial_t\de_n)| &\leq \frac{1}{64}\|\partial_t\de_n\|^2+C\Delta t\|\bu_{ht}\|^2_{L^2(t_{n-1},t_n;\varepsilon)},\\
			|c(\de_{n-1},\bu_h^n-\bu^n,\partial_t\de_n)+c(\de_{n-1},\bu^n,\partial_t\de_n)| &\leq \frac{1}{64}\|\partial_t\de_n\|^2+C \|\de_{n-1}\|_\varepsilon^2.
		\end{align}
		The last term in the right hand side of (\ref{ent3}) can be rewritten as
		\begin{align*}
			-c(\de_{n-1},\de_n,\partial_t\de_n)=-\frac{1}{\Delta t}c(\de_{n-1},\de_n,\de_n)+\frac{1}{\Delta t}c(\de_{n-1},\de_n,\de_{n-1}).
		\end{align*}
		From the positivity property (\ref{5.5}), we have $\frac{1}{\Delta t}c(\de_{n-1},\de_n,\de_n)\geq 0$. The estimate (\ref{tri2}) yields
		\begin{align}
			\frac{1}{\Delta t}|c(\de_{n-1},\de_n,\de_{n-1})|\leq \frac{C}{\Delta t}\|\de_{n-1}\|_\varepsilon^2\|\de_n\|_\varepsilon.
		\end{align}
		From Proposition 4.10 in \cite{GR09}, we can easily see that
		\begin{align*}
			|c^{\bu_h^n} &(\bU^{n-1},\bU^n,\partial_t\de_n)-c^{\bU^{n-1}}(\bU^{n-1},\bU^n,\partial_t\de_n)|\\
			&\leq \|\bu_h^n-\bU^{n-1}\|_{L^4(\Omega)}\|\bU^n-\bu^n\|_\varepsilon\|\partial_t\de_n\|_{L^4(\Omega)}\\
			& \leq \|\bu_h^n-\bU^{n-1}\|_\varepsilon\|\bu_h^n-\bu^n\|_\varepsilon \frac{1}{\displaystyle\min_{T\in \mathcal{T}_h}h_T^{1/2}}\|\partial_t\de_n\|+\|\bu_h^n-\bU^{n-1}\|^{1/2}\|\bu_h^n-\bU^{n-1}\|_\varepsilon^{1/2}\|\de_n\|_\varepsilon \|\partial_t\de_n\|^{1/2}\|\partial_t\de_n\|_\varepsilon^{1/2}.
		\end{align*}
		A use of the triangle inequality, Theorem \ref{semierror}, {\it a priori} estimate (\ref{uh1}), Lemma \ref{fullypriorisol} and the Cauchy-Schwarz inequality yield
		\begin{align}
			| &c^{\bu_h^n}(\bU^{n-1},\bU^n,\partial_t\de_n)-c^{\bU^{n-1}}(\bU^{n-1},\bU^n,\partial_t\de_n)|\nonumber\\
			\leq & \frac{1}{64}\|\partial_t\de_n\|^2+ C\Delta t\|\bu_{ht}\|_{L^2(t_{n-1},t_n;\varepsilon)}^2+C\|\de_{n-1}\|_\varepsilon^2 +\frac{C}{\Delta t} \|\bu_h^n-\bU^{n-1}\|_\varepsilon^{1/2}\|\de_n\|_\varepsilon(\|\de_n\|_\varepsilon^{1/2}+\|\de_{n-1}\|_\varepsilon^{1/2})\nonumber\\
			\leq & \frac{1}{64}\|\partial_t\de_n\|^2+ C\Delta t\|\bu_{ht}\|_{L^2(t_{n-1},t_n;\varepsilon)}^2+C\|\de_{n-1}\|_\varepsilon^2+ C\|\bu_{ht}\|_{L^2(t_{n-1},t_n;\varepsilon)}^2+\frac{C}{\Delta t}(\|\de_{n-1}\|_\varepsilon^2+\|\de_n\|_\varepsilon^2).
		\end{align}
		From (\ref{11.1}), we have 
		\begin{align}
			(\bu^n_{ht},\partial_t\de_n)-(\partial_t\bu_h^n,\partial_t\de_n)\leq C\Delta t^{1/2}\left(\int_{t_{n-1}}^{t_n}\|\bu_{htt}(t)\|_{-1,h}^2\,dt\right)^{1/2}\|\partial_t\de_n\|_\varepsilon\nonumber\\
			\leq  C\int_{t_{n-1}}^{t_n}\|\bu_{htt}(t)\|_{-1,h}^2\,dt+\frac{C}{\Delta t}\|\de_n-\de_{n-1}\|_\varepsilon^2.\label{ent9}
		\end{align}
		Since $a(\cdot,\cdot)$ is symmetric, one can obtain
		\begin{align}
			a(\de_n,\partial_t\de_n)= \frac{1}{2}\bigg(\frac{1}{\Delta t}a(\de_n,\de_n)-\frac{1}{\Delta t}a(\de_{n-1},\de_{n-1})+\Delta t a(\partial_t\de_n,\partial_t\de_n)\bigg).
		\end{align}
		Again,
		\begin{align}
			\sum_{n=1}^m  & e^{2\al t_n}\bigg((a+J_0)(\de_n,\de_n)-(a+J_0)(\de_{n-1},\de_{n-1})\bigg)\nonumber\\
			=&e^{2\al t_m}(a+J_0)(\de_m,\de_m)-\sum_{n=1}^{m-1} e^{2\al t_n}(e^{2\al \Delta t}-1)(a+J_0)(\de_n,\de_n).\label{ent10}
		\end{align}
		Combining (\ref{ent2})--(\ref{ent10}), multiply  (\ref{ent1}) by $\Delta t e^{2\al t_n}$, sum over $1\leq n\leq m\leq M$ and using Lemma \ref{coer}, we obtain
		\begin{align}
			\Delta t \sum_{n=1}^m & e^{2\al t_n}\|\partial_t\de_n\|^2+\mu K e^{2\al t_m}\|\de_m\|^2_\varepsilon\leq C\Delta t\sum_{n=1}^{m-1} e^{2\al t_n} (a+J_0)(\de_n,\de_n) +C\Delta t  \sum_{n=1}^m e^{2\al t_n}\|\de_n\|_\varepsilon^2\nonumber\\
			&+C\Delta t  \sum_{n=1}^m e^{2\al t_n}\|\de_{n-1}\|_\varepsilon^2 +\frac{C}{\Delta t}\Delta t  \sum_{n=1}^m e^{2\al t_n}\|\de_n\|_\varepsilon^2\|\de_{n-1}\|_\varepsilon^2+\frac{C}{\Delta t}\Delta t  \sum_{n=1}^m e^{2\al t_n}(\|\de_n\|_\varepsilon^2+\|\de_{n-1}\|_\varepsilon^2)\nonumber\\
			&+C{\Delta t}\sum_{n=1}^m e^{2\al t_n}\int_{t_{n-1}}^{t_n}\|\bu_{htt}(t)\|_{-1,h}^2\,dt+ C{\Delta t}\sum_{n=1}^m e^{2\al t_n}\int_{t_{n-1}}^{t_n}\|\bu_{ht}(t)\|_\varepsilon^2\,dt.\label{ent11}
		\end{align} 
		Using Lemmas \ref{cont}, \ref{uhpriori} and  \ref{fully}, from (\ref{ent11}), we obtain
		\begin{align*}
			\Delta t \sum_{n=1}^m  e^{2\al t_n}\|\partial_t\de_n\|^2+\mu K e^{2\al t_m}\|\de_m\|^2_\varepsilon\leq \frac{C}{\Delta t}\Delta t  \sum_{n=1}^m e^{2\al t_n}\|\de_n\|_\varepsilon^2\|\de_{n-1}\|_\varepsilon^2+C{\Delta t}e^{2\al t_{m+1}}.
		\end{align*} 
		Finally, a use of the discrete Gronwall's inequality and  Lemma \ref{fully} give us the desired estimate. This completes the proof. 
	\end{proof}
	
	\begin{lemma}\label{neget}
		The error $\de_n=\bU^n-\bu_h^n$, satisfies
		\begin{align*}
			\|\partial_t\de_n\|_{-1,h}\leq C\Delta t^{1/2}.
		\end{align*}
	\end{lemma}
	\begin{proof}
		The non-linear terms in the error equation (\ref{11.0}) can be rewritten as
		\begin{align}\label{neget2}
			c^{\bu_h^n}(\bu_h^n,\bu_h^n,\bv_h)-c^{\bU^{n-1}}(\bU^{n-1},\bU^n,\bv_h)=&-c^{\bu_h^n}(\bu_h^{n-1},\de_n,\bv_h) +c^{\bu_h^n}(\bu_h^n-\bu_h^{n-1},\bu_h^n-\bu^n,\bv_h)\nonumber\\
			&+c^{\bu_h^n}(\bu_h^n-\bu_h^{n-1},\bu^n,\bv_h) -c^{\bu_h^n}(\de_{n-1},\de_n,\bv_h)-c^{\bu_h^n}(\de_{n-1},\bu_h^n,\bv_h)\nonumber\\
			& +c^{\bu_h^n}(\bU^{n-1},\bU^n,\bv_h)-c^{\bU^{n-1}}(\bU^{n-1},\bU^n,\bv_h).
		\end{align}
		Using estimate (\ref{tri2}), we obtain
		\begin{align}\label{neget3}
			|c^{\bu_h^n}(\bu_h^{n-1},\de_n,\bv_h)+c^{\bu_h^n}(\de_{n-1},\de_n,\bv_h)+c^{\bu_h^n}(\de_{n-1},\bu_h^n,\bv_h)|\leq C\|\bu_h^{n-1}\|_\varepsilon\|\de_n\|_\varepsilon\|\bv_h\|_\varepsilon+C\|\de_{n-1}\|_\varepsilon\|\de_n\|_\varepsilon\|\bv_h\|_\varepsilon\nonumber\\
			+C\|\de_{n-1}\|_\varepsilon\|\bu_h^n\|_\varepsilon\|\bv_h\|_\varepsilon.
		\end{align}
		Again, by using  Theorem \ref{semierror}, estimate (\ref{tri1}) and  Sobolev inequality, one can obtain the bounds similar to Lemma \ref{fully} as follows:
		\begin{align}\label{neget4}
			|c^{\bu_h^n}(\bu_h^n-\bu_h^{n-1},\bu_h^n-\bu^n,\bv_h)+c^{\bu_h^n}(\bu_h^n-\bu_h^{n-1},\bu^n,\bv_h)|\leq C\|\bu_h^n-\bu_h^{n-1}\|\|\bv_h\|_\varepsilon.
		\end{align}
		Furthermore, similar to Lemma \ref{L2L2delt} and using Theorem \ref{semierror}, we find that
		\begin{align}\label{neget5}
			|c^{\bu_h^n}(\bU^{n-1} &,\bU^n,\bv_h)-c^{\bU^{n-1}}(\bU^{n-1},\bU^n,\bv_h)|\nonumber\\
			\leq &\|\bu_h^n-\bU^{n-1}\|_{L^4(\Omega)}\|\bU^n-\bu^n\|_\varepsilon\|\bv_h\|_{L^4(\Omega)}\nonumber\\
			\leq & \|\bu_h^n-\bU^{n-1}\|_{L^4(\Omega)}\|\de_n\|_\varepsilon\|\bv_h\|_\varepsilon+ \|\bu_h^n-\bU^{n-1}\|_{L^4(\Omega)}\|\bu^n-\bu_h^n\|_\varepsilon\|\bv_h\|_\varepsilon\nonumber\\
			\leq & \|\bu_h^n-\bu_h^{n-1}\|_\varepsilon\|\de_n\|_\varepsilon\|\bv_h\|_\varepsilon+\|\de_{n-1}\|_\varepsilon\|\de_n\|_\varepsilon\|\bv_h\|_\varepsilon+C\|\bu_h^n-\bu_h^{n-1}\|\|\bv_h\|_\varepsilon+C\|\de_{n-1}\|\|\bv_h\|_\varepsilon.
		\end{align}
		Now, Lemma \ref{cont} yield
		\begin{align}\label{neget6}
			|a(\de_n,\bv_h)+J_0(\de_n,\bv_h)|\leq C\|\de_n\|_\varepsilon\|\bv_h\|_\varepsilon.
		\end{align}
		Applying (\ref{11.1}), Cauchy-Schwarz's inequality, Young's inequality and estimate (\ref{uh2}), we arrive at 
		\begin{align}
			(\bu^n_{ht},\de_n)-(\partial_t\bu_h^n,\bv_h)& \leq C\Delta t^{1/2}\left(\int_{t_{n-1}}^{t_n}\|\bu_{htt}(t)\|^2_{-1,h}\,dt\right)^{1/2}\|\bv_h\|_\varepsilon\leq C\Delta t^{1/2}\|\bv_h\|_\varepsilon.\label{neget7}
		\end{align}
		Combining all the bounds (\ref{neget2})--(\ref{neget7}) in (\ref{11.0}), using the definition of $\|\cdot\|_{-1,h}$  and finally using Lemmas \ref{uhpriori} and \ref{L2L2delt} , we obtain our desired result. This completes the rest of the proof.
	\end{proof}
	\begin{lemma}\label{fullypressure}
		There exists a positive constant $C$, independent of $\Delta t$, such that for $n=1,2,\cdots,M$
		\begin{align*}
			\|P^n-p_h^n\|\leq C\Delta t^{1/2}.
		\end{align*}
	\end{lemma}
	\begin{proof}
		Consider (\ref{8.8}) at $t=t_n$ and subtract it from (\ref{fullyxh1}) to obtain
		\begin{align}
			b(\bv_h,P^n-p_h^n)= &(\bu_{ht}^n,\bv_h) -(\partial_t\bu_h^n,\bv_h)-(\partial_t\de_n,\bv_h)-\mu(  a(\de_n,\bv_h)+J_0(\de_n,\bv_h))\nonumber\\
			&+c^{\bu_h^n}(\bu_h^n,\bu_h^n,\bv_h)-c^{\bU^{n-1}}(\bU^{n-1},\bU^n,\bv_h),\quad \forall \bv_h\in \bX_h.\label{fullyerrpn1}
		\end{align}
		Using Lemma \ref{inf-sup} and bound the terms on the right hand side of (\ref{fullyerrpn1}) following the steps involved in the proof of Lemma \ref{neget}, we arrive at
		\begin{align}
			\|P^n-p_h^n\|\leq C \|\partial_t\de_n\|_{-1,h}+C\|\de_n\|_\varepsilon+C\Delta t^{1/2}.\label{fullyerrpn2}
		\end{align}
		Finally, a use of the  Lemmas \ref{neget} and \ref{L2L2delt}, we conclude the rest of the proof. 
	\end{proof}
	\noindent
	A combination of Lemma \ref{fullypressure} and Theorem \ref{semipressure} lead to the following fully discrete error estimate for pressure.
	\begin{theorem}
		Under the assumptions of Theorem \ref{semipressure} and Lemma \ref{fullypressure}, the following hold true:
		\begin{align*}
			\|p(t_n)-P^n\|\leq C(ht_n^{-1/2}+\Delta t^{1/2}),
		\end{align*}
		where $C>0$ is a constant, independent of $h$ and $\Delta t$.
	\end{theorem}
	\section{Numerical experiments}\label{s7}
\begin{figure}[h!]
		\centering
		\includegraphics[scale=0.5]{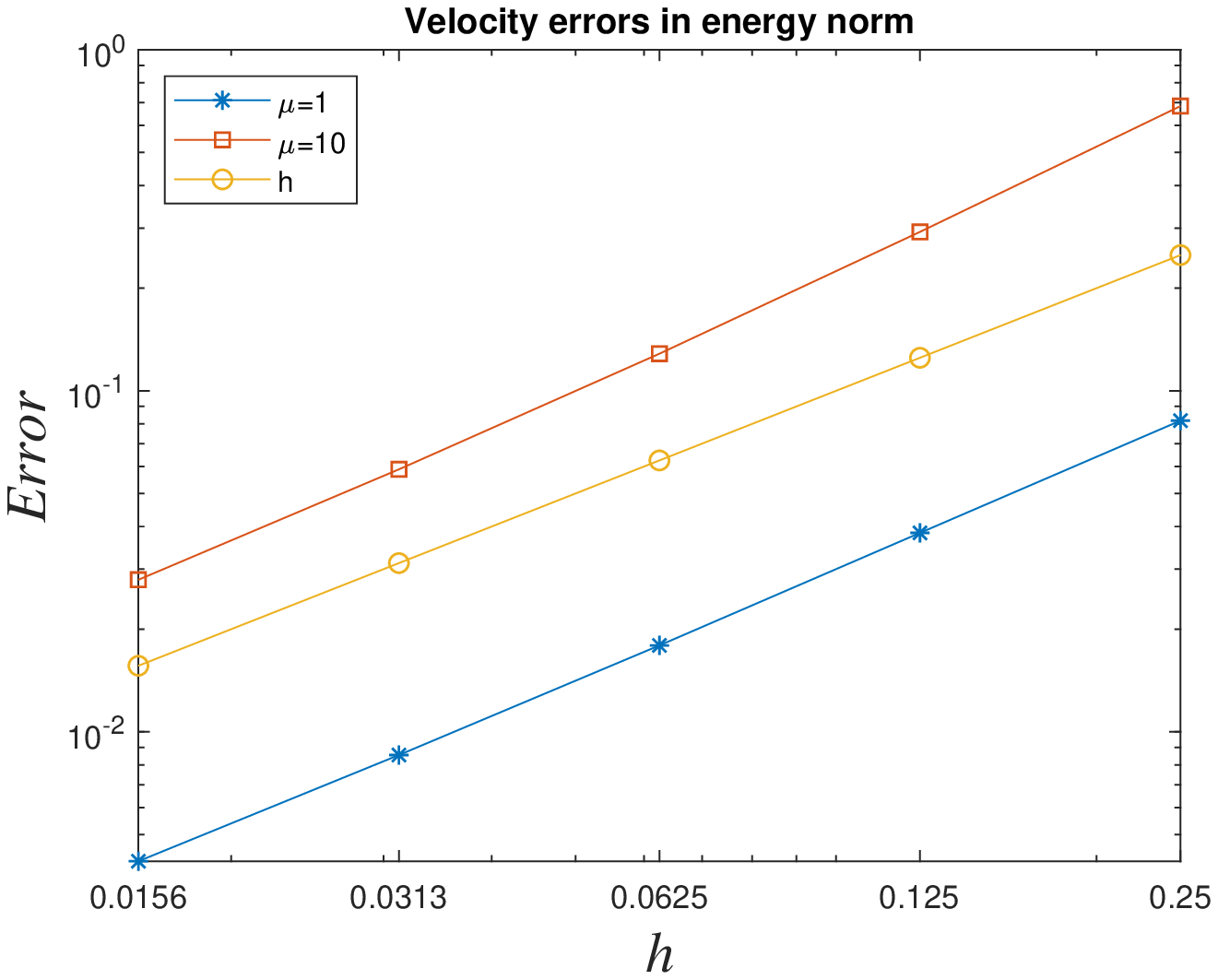}
		\includegraphics[scale=0.5]{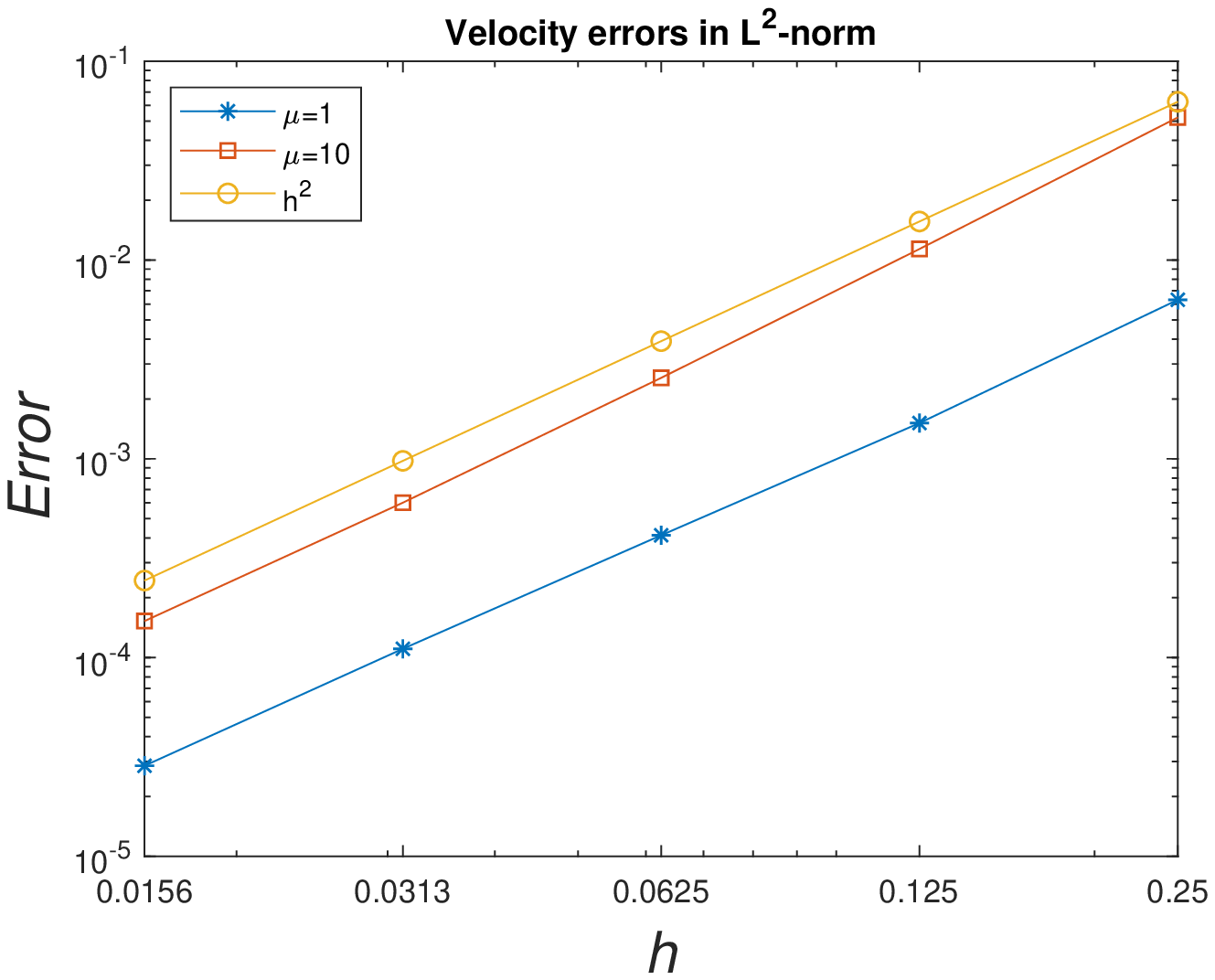}
		\includegraphics[scale=0.5]{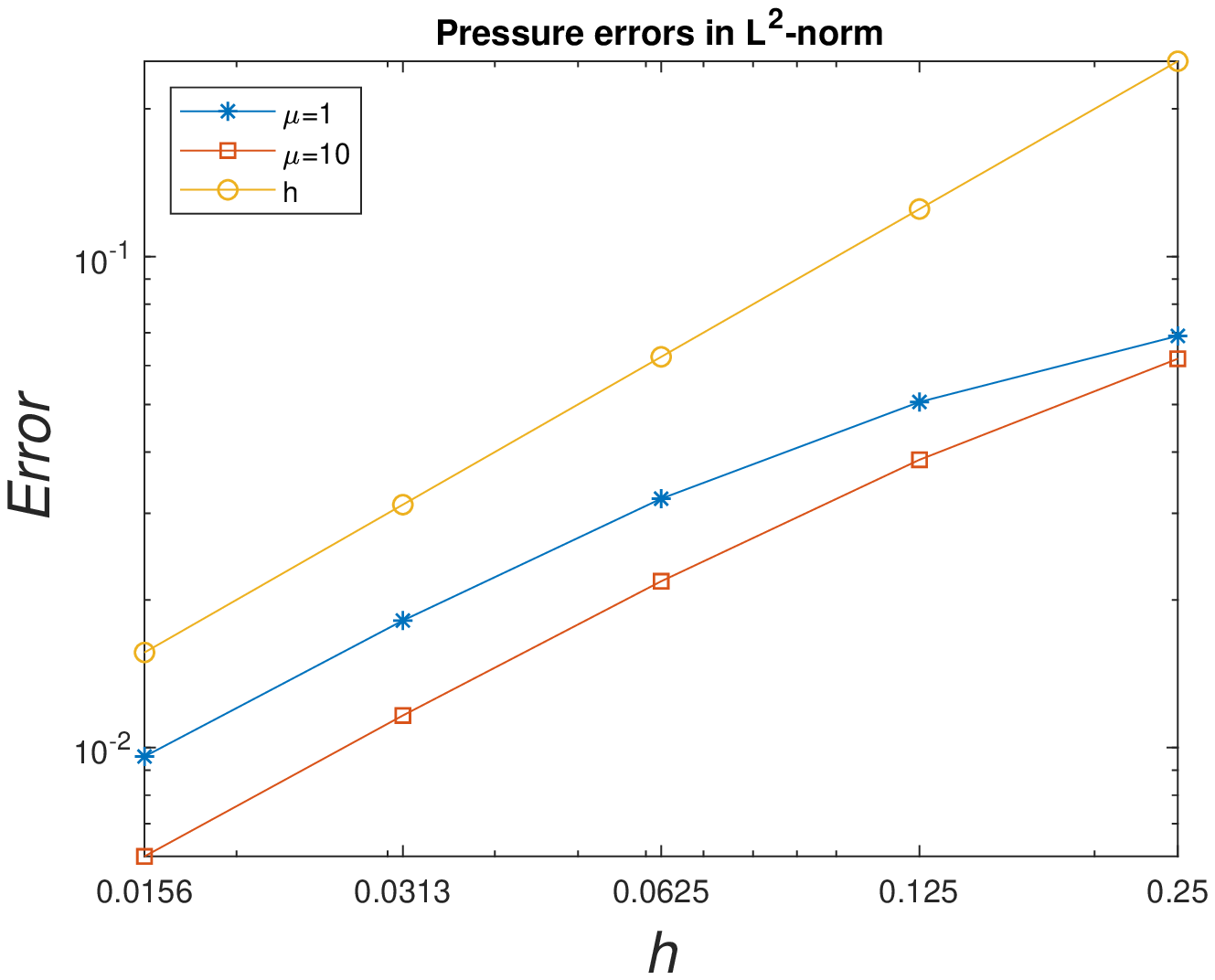}
		\caption{Velocity and pressure errors for $P_1-P_0$  with $\mu=1$ and $1/10$ for Example \ref{ex1}.}
		\label{fig1}
	\end{figure}

	\se
	In this section, we present numerical experiments to support the theoretical results.
	For space discretization, $P_1-P_0$ and $P_1-P_1$  mixed finite element spaces are used and for the time discretization, a first order accurate backward Euler method is applied. The spatial domain $\Omega$ is chosen as $ (0, 1)\times(0, 1)$ and the time interval is chosen as $[0,1]$ with final time $T=1$.
	\begin{example}\label{ex1}
		Consider the transient Navier-Stokes equations (\ref{8.1})-(\ref{8.4}) with exact solution as
		\[\bu = (2e^t x^2 (x-1)^2 y (y-1) (2y-1),-2e^t x(x-1)(2x-1)y^2(y-1)^2),\quad p=2e^t(x-y).\]
	\end{example}
	Tables 1 and 2 show the errors and the convergence rates for the mixed finite element space $P_1-P_0$  with viscosity $\mu = 1$ and $\mu = 1/10$, respectively. And Figure \ref{fig1} represents the velocity and pressure errors for $P_1-P_0$ element with $\mu=1$ and $10$.
	In Table 3, we choose $P_1-P_1$ mixed finite element space with $\mu = 1$. We choose a constant penalty parameter $\sigma_e = 10$ for Tables 1-3. 
	In Table 4, we have represented the errors and the convergence rates for the backward Euler method applied to continuous Galerkin finite element method with $\mu=1$.
	The numerical results represented in Tables 1-3 validate our theoretical findings in Theorem \ref{velofinal}. Further, the results in Tables 3 and 4 represent that the discontinuous Galerkin finite element method works well for equal order element $P_1-P_1$ whereas continuous Galerkin finite element method fails to approximate the exact solution \cite{GR79}. 
	{\small
		\renewcommand{\arraystretch}{1.3}
		\begin{table}[ht!]
			\centering
			\caption{Numerical errors and convergence rates, for $P_1$--$P_0$ ($\mu$=1, $\sigma_e=10$) for Example \ref{ex1} with $\Delta t=\mathcal{O}(h^2)$.}
			\vspace{.4cm}
			\begin{tabular}{|c|c c| c c| c c|}
				
				\hline
				
				$h$ & $\|\bu(T)-\bU^M\|_\varepsilon$ & Rate & $\|\bu(T)-\bU^M\|_{L^2(\Omega)}$ & Rate & $\|p(T)-P^M\|_{L^2(\Omega)}$ & Rate\\
				
				\hline  
				1/4  & $8.1759\times 10^{-2}$  &  		  & $6.3073\times 10^{-3}$ & 		 & $6.8941\times 10^{-2}$ & 		\\
				1/8  & $3.8398\times 10^{-2}$  &  1.0932  & $1.5131\times 10^{-3}$ & 2.0594 & $5.0580\times 10^{-2}$ & 0.4468\\
				1/16 & $1.7926\times 10^{-2}$  &  1.0960  & $4.1238\times 10^{-4}$ & 1.8754 & $3.2138\times 10^{-2}$ & 0.6542 \\
				1/32 & $8.5526\times 10^{-3}$  &  1.0676  & $1.1070\times 10^{-4}$ & 1.8973 & $1.8147\times 10^{-2}$ & 0.8244 \\
				1/64 & $4.1749\times 10^{-3}$  &  1.0346  & $2.8571\times 10^{-5}$ & 1.9540 & $9.5931\times 10^{-3}$ & 0.9197 \\
				\hline
				
			\end{tabular}
			\vspace{.1cm}
			
		\end{table}
	}

	{\small
		\renewcommand{\arraystretch}{1.3}
		\begin{table}[ht!]
			\centering
			\caption{Numerical errors and convergence rates, for $P_1$--$P_0$ ($\mu=1/10$, $\sigma_e$=10) for Example \ref{ex1} with $\Delta t=\mathcal{O}(h^2)$.}
			\vspace{.4cm}
			\begin{tabular}{|c|c c| c c| c c|}
				
				\hline
				
				$h$ & $\|\bu(T)-\bU^M\|_\varepsilon$ & Rate & $\|\bu(T)-\bU^M\|_{L^2(\Omega)}$ & Rate & $\|p(T)-P^M\|_{L^2(\Omega)}$ & Rate\\
				
				\hline 
				1/4 & $0.68355$   &  	  & $5.2164\times 10^{-2}$ & 	 & $6.1942\times 10^{-2}$ &	 \\
				1/8  & $0.29241$  & 1.2250  & $1.1363\times 10^{-2}$ & 2.1986 & $3.8573\times 10^{-2}$ & 0.6833\\
				1/16 & $0.12847$  & 1.1865  & $2.5536\times 10^{-3}$ & 2.1537 & $2.1817\times 10^{-2}$ & 0.8220\\
				1/32 & $0.05884$  & 1.1264  & $6.0150\times 10^{-4}$ & 2.0859 & $1.1627\times 10^{-2}$ & 0.9079 \\
				1/64 & $0.02792$  & 1.0753  & $1.5268\times 10^{-4}$ & 1.9780 & $6.0054\times 10^{-3}$ & 0.9531\\
				\hline
				
			\end{tabular}
			\vspace{.1cm}
			
		\end{table}
	}
	
	{\small
		\renewcommand{\arraystretch}{1.3}
		\begin{table}[ht!]
			\centering
			\caption{Numerical errors and convergence rates, for $P_1$--$P_1$ ($\mu=1$, $\sigma_e$=10) for Example \ref{ex1} with $\Delta t=\mathcal{O}(h^2)$.}
			\vspace{.4cm}
			\begin{tabular}{|c|c c| c c| c c|}
				
				\hline
				
				$h$ & $\|\bu(T)-\bU^M\|_\varepsilon$ & Rate & $\|\bu(T)-\bU^M\|_{L^2(\Omega)}$ & Rate & $\|p(T)-P^M\|_{L^2(\Omega)}$ & Rate\\
				
				\hline   
				1/4  & $0.04558$ &  		 & $0.003248$ & 		 & $0.1609$ &		\\
				1/8  & $0.02739$  &  0.7349  & $0.001614$ & 1.0085 & $0.1140$ & 0.4969\\
				1/16 & $0.01263$  &  1.1158 & $0.000589$ & 1.4545 & $0.0636$ & 0.8423\\
				1/32 & $0.00575$  &  1.1337 & $0.000169$ & 1.7934 & $0.0334$ & 0.9262 \\
				1/64 & $0.00275$ &  1.0651  & $0.0000486$ & 1.9215 & $0.0171$ & 0.9672\\
				\hline
				
			\end{tabular}
			\vspace{.1cm}
			
		\end{table}
	}

	{\small
		\renewcommand{\arraystretch}{1.3}
		\begin{table}[ht!]
			\centering
			\caption{Numerical errors and convergence rates for continuous finite element method using $P_1$--$P_1$ ($\mu=1$) for Example \ref{ex1} with $\Delta t=\mathcal{O}(h^2)$.}
			\vspace{.4cm}
			\begin{tabular}{|c|c c| c c| c c|}
				
				\hline
				
				$h$ & $\|\bu(T)-\bU^M\|_{H^1(\Omega)}$ & Rate & $\|\bu(T)-\bU^M\|_{L^2(\Omega)}$ & Rate & $\|p(T)-P^M\|_{L^2(\Omega)}$ & Rate\\
				
				\hline   
				1/4 & $9.1515\times 10^{-2}$ & 		 & $9.2951\times 10^{-3}$ & 		& $1.07149$ &	\\
				1/8  & $4.1709\times 10^{-2}$  &  1.1336 & $2.5728\times 10^{-3}$ & 1.8531 & $0.59502$ & 0.8486\\
				1/16 & $1.9316\times 10^{-2}$  &  1.1105 & $6.2603\times 10^{-4}$ & 2.0390 & $0.36037$ & 0.7234\\
				1/32 & $9.3341\times 10^{-3}$  &  1.0492 & $1.5252\times 10^{-4}$ & 2.0277 & $0.26159$ & 0.4621 \\
				1/64 & $4.5978\times 10^{-3}$  & 1.0215  & $3.7996\times 10^{-5}$ & 2.0145 & $0.22689$ &  0.2053\\
				\hline
				
			\end{tabular}
			\vspace{.1cm}
			
		\end{table}
	}

\begin{figure}[h!]
		\centering
		\includegraphics[scale=0.5]{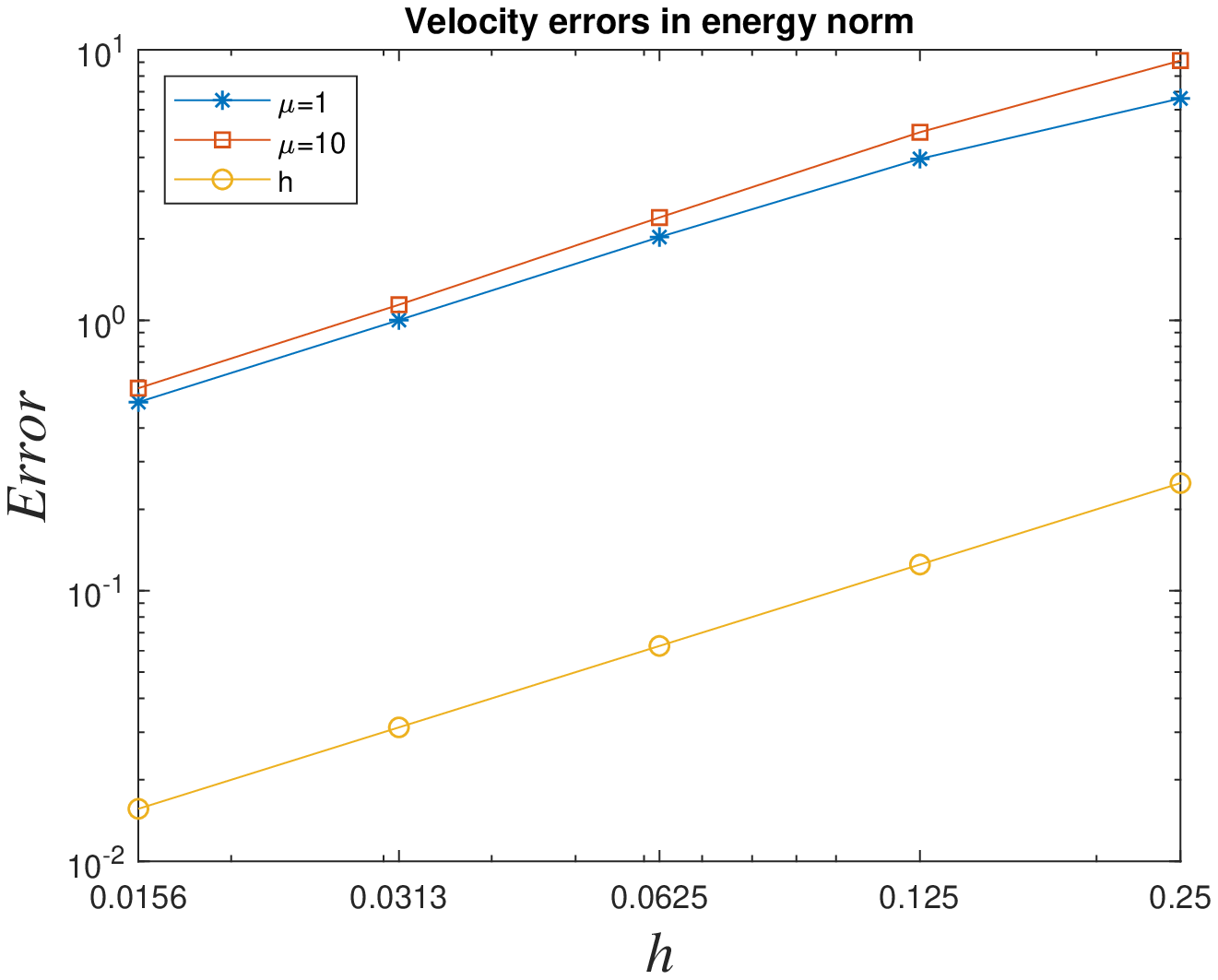}
		\includegraphics[scale=0.5]{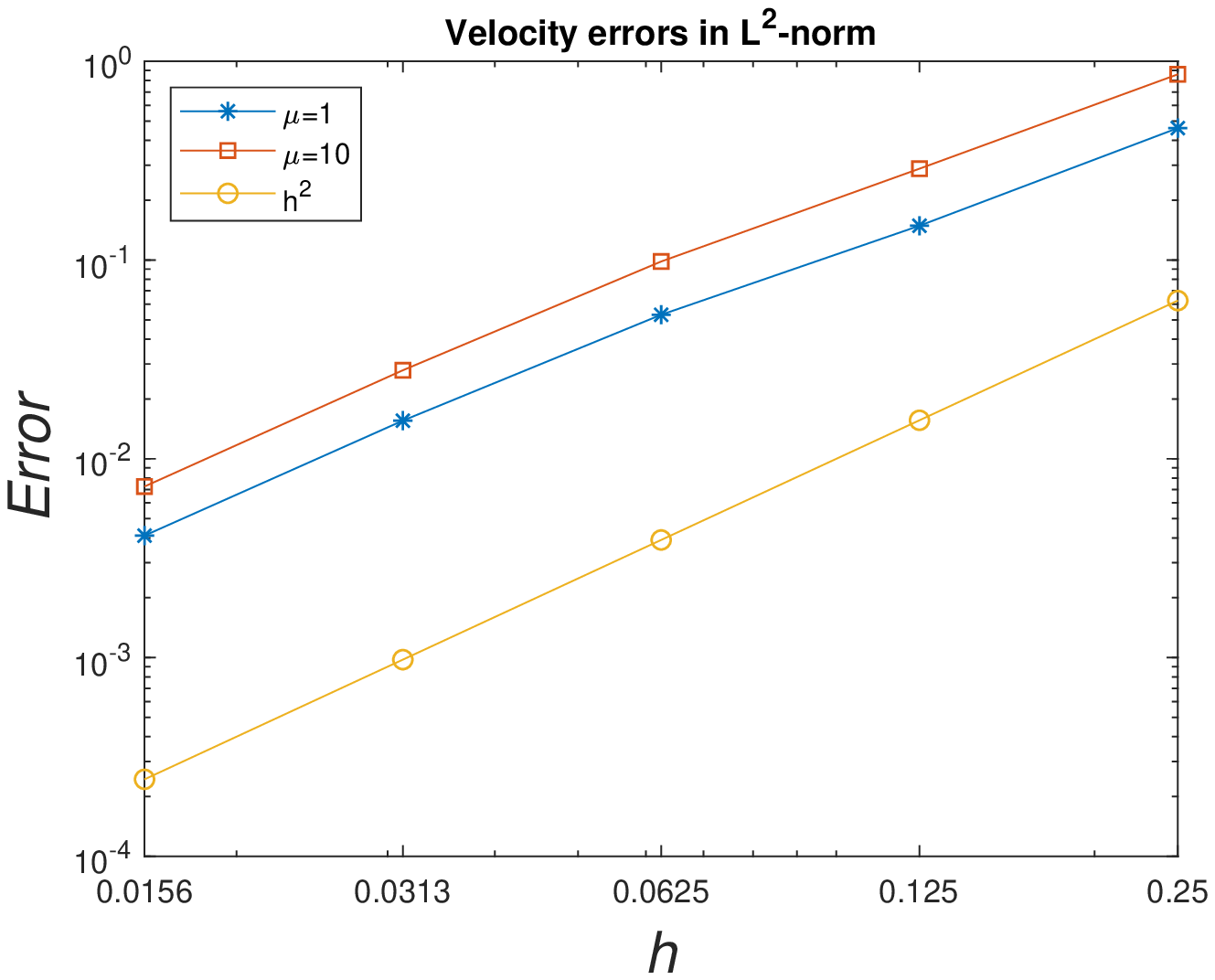}
		\includegraphics[scale=0.5]{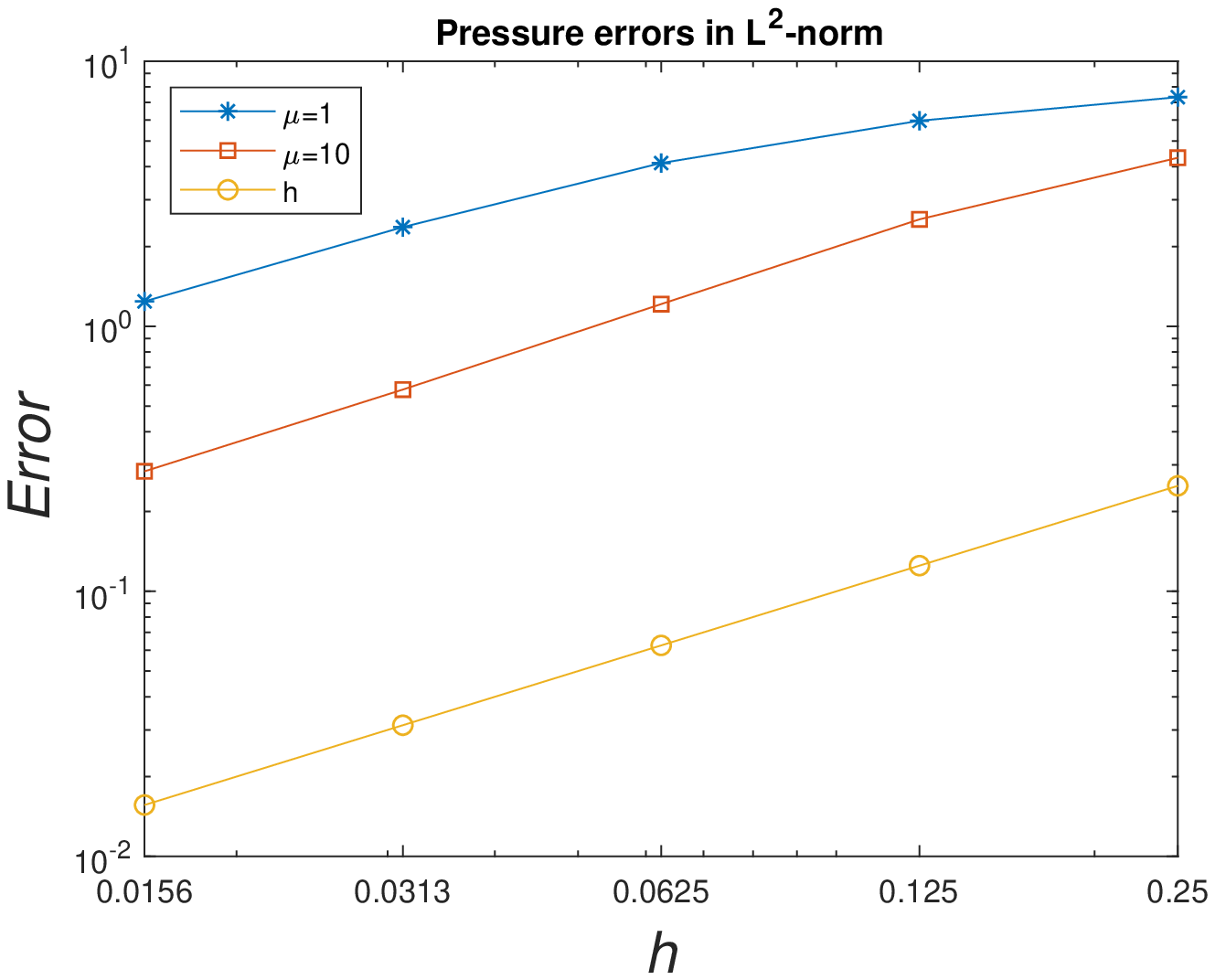}
		\caption{Velocity and pressure errors for $P_1-P_0$  with $\mu=1$ and $1/10$ for Example \ref{ex2}.}
		\label{fig2}
	\end{figure}

	\begin{example}\label{ex2}
		In this example, we choose the right hand side function f in such a way
		that the exact solution is:
		\[\bu = (e^t(-cos(2\pi x)*sin(2\pi y) + sin(2\pi y)),e^t(sin(2\pi x)cos(2\pi y) - sin(2\pi x)),\quad p=e^t (2\pi (cos(2\pi y)- cos(2\pi x))).\]
	\end{example}
	In Tables 5 and 6, we have shown the error and convergence rates for the mixed finite element space $P_1-P_0$  with viscosity $\mu = 1$ and $\mu = 1/10$ respectively.  And Figure \ref{fig2} represents the velocity and pressure errors for $P_1-P_0$ element with $\mu=1$ and $10$.
	For the cases $\mu = 1$ and $\mu = 1/10$, we have choosen  $\sigma_e=10$ and  $\sigma_e=20$ respectively. Table 7 depicts the result for the mixed space $P_1-P_1$ with $\mu = 1$ and $\sigma_e=10$.
	For continuous finite element method, the errors and the convergence rates for the case $P_1-P_1$ with $\mu=1$ have shown in Table 8 for Example \ref{ex2}.
	{\small
		\renewcommand{\arraystretch}{1.3}
		\begin{table}[ht!]
			\centering
			\caption{Numerical errors and convergence rates, for $P_1$--$P_0$ ($\mu=1$, $\sigma_e$=10) for Example \ref{ex2} with $\Delta t=\mathcal{O}(h^2)$.}
			\vspace{.4cm}
			\begin{tabular}{|c|c c| c c| c c|}
				
				\hline

				$h$ & $\|\bu(T)-\bU^M\|_\varepsilon$ & Rate & $\|\bu(T)-\bU^M\|_{L^2(\Omega)}$ & Rate & $\|p(T)-P^M\|_{L^2(\Omega)}$ & Rate\\
				
				\hline   
				1/4  & $6.6089$  &  	  & $0.4609$ & 		 & $7.3185$ &	\\
				1/8  & $3.9510$  &  0.7421  & $0.1487$ & 1.6317 & $5.9624$ & 0.2956\\
				1/16 & $2.0300$  &  0.9607 & $0.0530$ & 1.4876 & $	4.1285$ & 0.5302\\
				1/32 & $1.0011$  &  1.0197 & $0.0155$ & 1.7693 & $2.3667$ & 0.8027 \\
				1/64 & $0.4979$  &  1.0075  & $0.0041$ & 1.9213 & $1.2414$ & 0.9308\\
				\hline
			\end{tabular}
			\vspace{.1cm}
			
		\end{table}
	}
	
	{\small
		\renewcommand{\arraystretch}{1.3}
		\begin{table}[ht!]
			\centering
			\caption{Numerical errors and convergence rates, for $P_1$--$P_0$ ($\mu=1/10$ and $\sigma_e$=20) for Example \ref{ex2} with $\Delta t=\mathcal{O}(h^2)$.}
			\vspace{.4cm}
			\begin{tabular}{|c|c c| c c| c c|}
				
				\hline
				
				$h$ & $\|\bu(T)-\bU^M\|_\varepsilon$ & Rate & $\|\bu(T)-\bU^M\|_{L^2(\Omega)}$ & Rate & $\|p(T)-P^M\|_{L^2(\Omega)}$ & Rate\\
				
				\hline   
				1/4  & $9.11955$  &         & $0.85917$ &        & $4.32282$ &            \\
				1/8  & $4.95544$  &  0.8799  & $0.28801$ & 1.5768 & $2.53199$ & 0.7716\\
				1/16 & $2.39775$  &  1.0473 & $0.09829$ & 1.5509 & $1.21246$ & 1.0623\\
				1/32 & $1.14083$  &  1.0715 & $0.02788$ & 1.8174 & $0.57703$ & 1.0712\\
				1/64 & $0.56084$  & 1.0244  & $0.00725$ &1.9425 & $0.28362$ & 1.0246\\
				\hline
			\end{tabular}
			\vspace{.1cm}
			
		\end{table}
	}

	{\small
		\renewcommand{\arraystretch}{1.3}
		\begin{table}[ht!]
			\centering
			\caption{Numerical errors and convergence rates, for $P_1$--$P_1$ ($\mu=1$, $\sigma_e$=10) for Example \ref{ex2} with $\Delta t=\mathcal{O}(h^2)$.}
			\vspace{.4cm}
			\begin{tabular}{|c|c c| c c| c c|}
				
				\hline
				
				$h$ & $\|\bu(T)-\bU^M\|_\varepsilon$ & Rate & $\|\bu(T)-\bU^M\|_{L^2(\Omega)}$ & Rate & $\|p(T)-P^M\|_{L^2(\Omega)}$ & Rate\\
				
				\hline   
				1/4  & $7.04721$  &  	     & $0.46170$ & 		 & $27.59425$ &   	\\
				1/8  & $4.47385$  &  0.6555 & $0.27692$ & 0.7374 & $19.82990$ & 0.4766\\
				1/16 & $2.02634$  &  1.1426 & $0.10272$ & 1.4306 & $10.55237$ & 0.9101\\
				1/32 & $0.90844$  &  1.1574 & $0.02933$ & 1.8082 & $5.42310$  & 0.9603 \\
				1/64 & $0.43356$  &  1.0671 & $0.00765$ & 1.9386 & $2.73855$  & 0.9857\\
				\hline
			\end{tabular}
			\vspace{.1cm}
			
		\end{table}
	}

	{\small
		\renewcommand{\arraystretch}{1.3}
		\begin{table}[ht!]
			\centering
			\caption{Numerical errors and convergence rates, for continuous finite element method using $P_1$--$P_1$ ($\mu=1$) for Example \ref{ex2} with $\Delta t=\mathcal{O}(h^2)$.}
			\vspace{.4cm}
			\begin{tabular}{|c|c c| c c| c c|}
				
				\hline
				
				$h$ & $\|\bu(T)-\bU^M\|_{H^1(\Omega)}$ & Rate & $\|\bu(T)-\bU^M\|_{L^2(\Omega)}$ & Rate & $\|p(T)-P^M\|_{L^2(\Omega)}$ & Rate\\
				
				\hline   
				1/4  & $14.49369$  &  		& $1.43048$ &  		& $68.66457$ &		\\
				1/8  & $4.48685$  &  1.6916 & $0.29520$ & 2.2767 & $30.17823$ & 1.1860\\
				1/16 & $1.91850$  &   1.2257 & $ 0.06874$ &  2.1024 & $23.51860$ & 0.3597\\
				1/32 & $ 0.90913$  &  1.0774 & $0.01672$ & 2.0392 & $21.55581$  & 0.1257 \\
				1/64 & $0.44585$  &  1.0279 & $ 0.00413$ & 2.0162 & $ 20.82569$  &  0.0497\\
				\hline
			\end{tabular}
			\vspace{.1cm}
		\end{table}
	}

	\pagebreak
	
\begin{figure}[h!]
		\centering
		\includegraphics[scale=0.5]{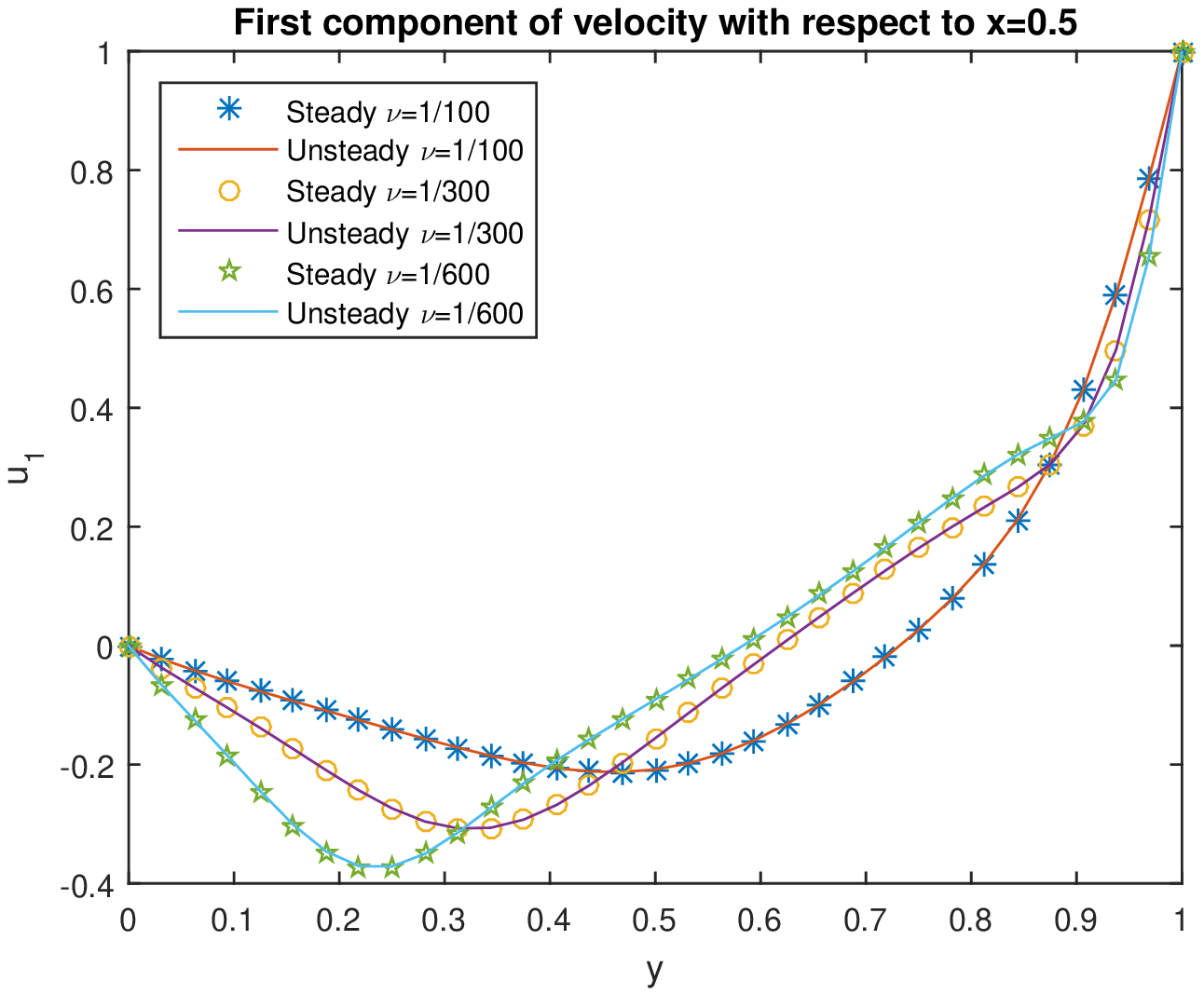}
		\includegraphics[scale=0.5]{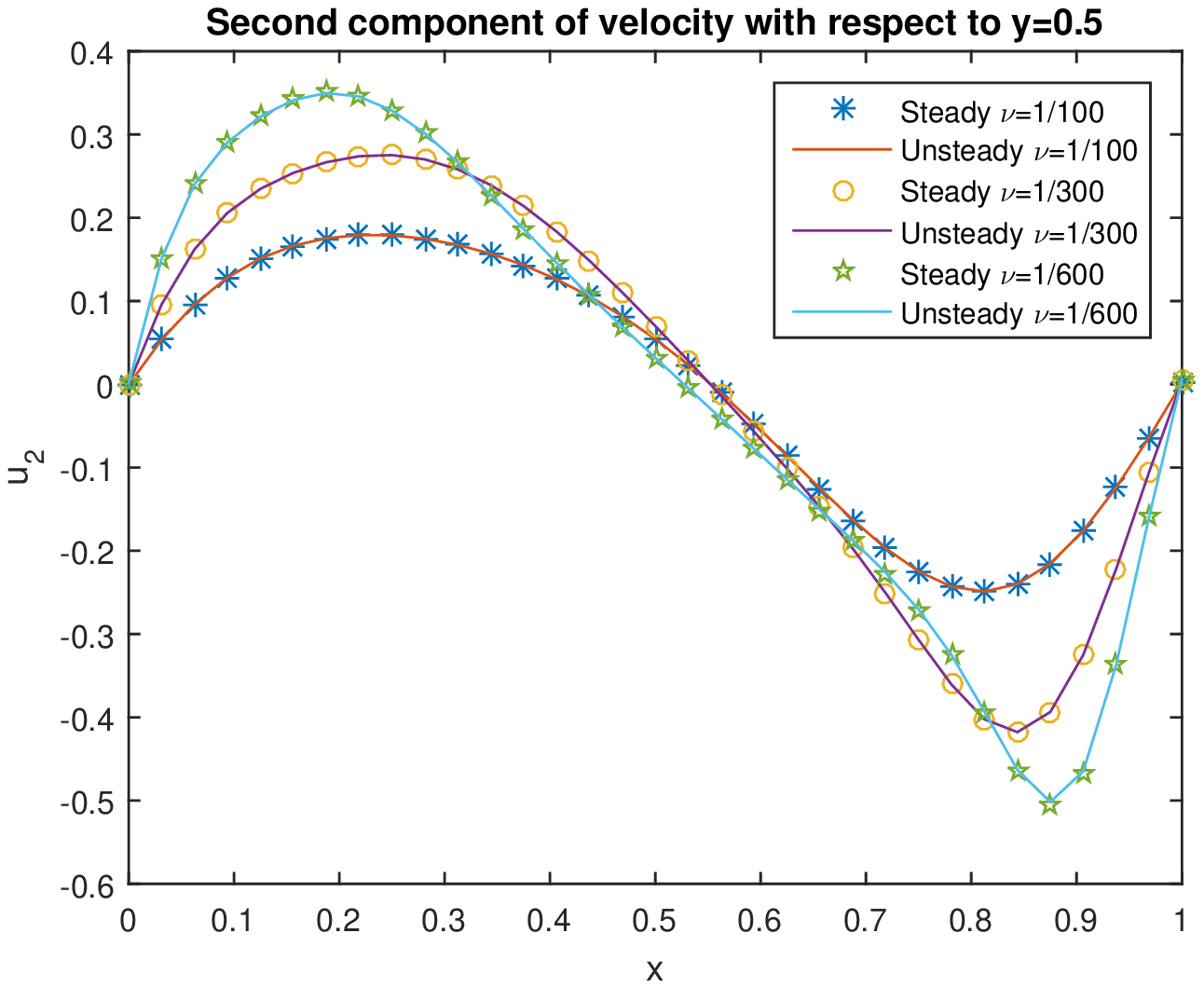}
		\caption{Velocity components for Example \ref{ex3}.}
		\label{fig3}
	\end{figure}

\begin{figure}[h!]
		\centering
		\includegraphics[scale=0.5]{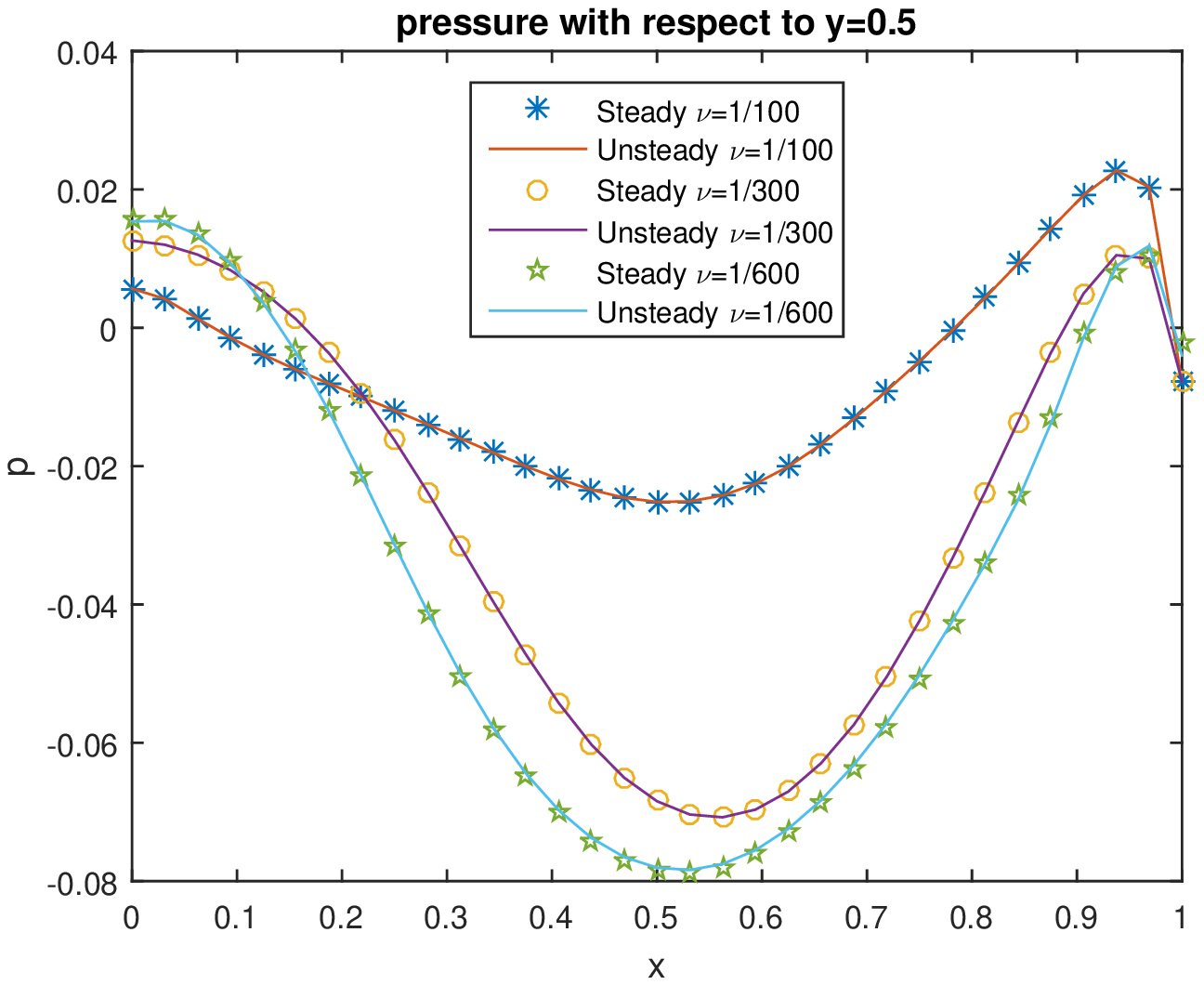}
		\includegraphics[scale=0.5]{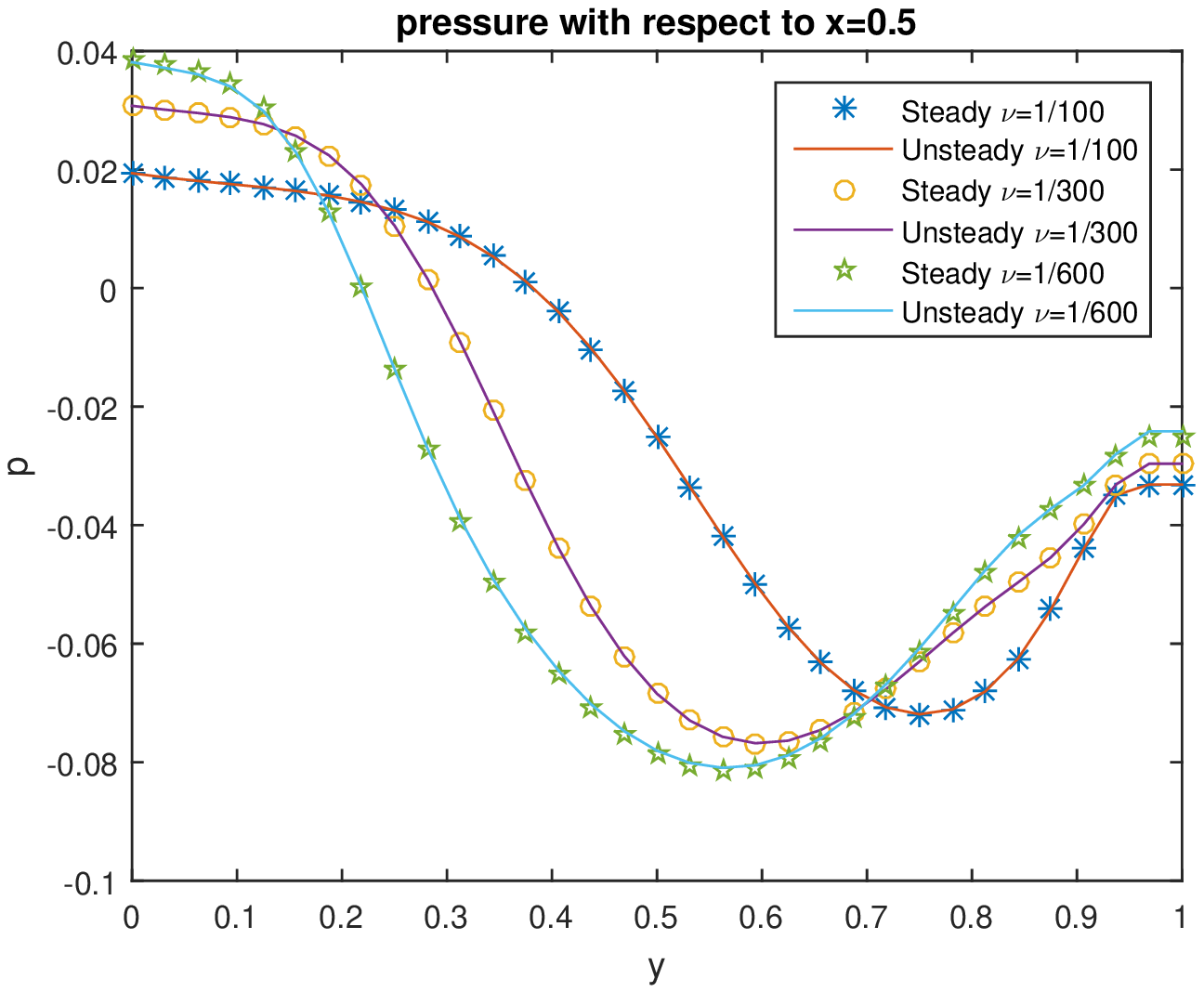}
		\caption{Pressure for Example \ref{ex3}.}
		\label{fig4}
	\end{figure}

	\begin{example}\label{ex3}
	(2D Lid Driven Cavity Flow Benchmark Problem). In this example, we consider a benchmark problem related to a lid driven cavity flow on a unit square with zero body forces. Further, no slip boundary conditions are considered everywhere except non zero velocity $(u_1,u_2)=(1,0)$ on the upper part of boundary, that is, the lid of the cavity is moving horizontally with a prescribed velocity. For numerical experiments, we have chosen lines $(0.5,y$) and $(x,0.5)$. We choose $P_1-P_0$ mixed finite element space for the space discretization. In Figure \ref{fig3}, we have
	presented the comparison between unstedy  backward Euler and steady state velocities, whereas Figure \ref{fig4} depicts the comparison of backward Euler and steady state pressure for different values of viscosity $\mu = \{1/100,1/300,1/600\}$, final time $t = 75$, $h = 1/32$ and $\Delta t = \mathcal{O}(h^2)$ with $\sigma_e=40$.  From the graphs, it is observed that the time dependent Navier-Stokes solution converges to the steady state solution for large time which validates the theoretical results.
	\end{example}
	\pagebreak
	\section*{References}
	\begingroup
	\renewcommand{\section}[2]{}%
	
\end{document}